\documentclass{article} 

\usepackage{authblk}

\usepackage[utf8]{inputenc} 
\usepackage[T1]{fontenc}    
 \usepackage{hyperref}  
\usepackage{hypcap}   
\usepackage{bookmark} 
\usepackage{url}            
\usepackage{booktabs}       
\usepackage{nicefrac}       
\usepackage{microtype}      
\usepackage{lipsum}
\usepackage{graphicx}
\graphicspath{ {./images/} }

\usepackage[a4paper, total={6in, 8in}]{geometry}

\usepackage[english]{babel}
\usepackage{tikz}
\usetikzlibrary{shapes.arrows}
\usepackage{xcolor}
\usepackage{pgfplots}
\pgfplotsset{compat=1.8}
\usepgfplotslibrary{statistics}

\usepackage{array}
\usepackage{longtable}
\usepackage{xurl}

\usepackage{amssymb,subfigure,icomma,xcolor}

\usepackage{amsmath}
\usepackage{amsthm}

\usepackage{multirow}
\usepackage{pgfplots}

\usepackage{multicol}

\providecommand{\keywords}[1]{\textbf{\textit{Keywords:}} #1}

\begin{document}

\title{Maximum flow-based formulation for the optimal location of electric vehicle charging stations}

\author[1]{Pierre-Luc Parent}

\author[1]{Margarida Carvalho}

\author[2]{Miguel F. Anjos}

\author[3]{Ribal Atallah}


\affil[1]{CIRRELT and Département d’informatique et de recherche opérationnelle, Université de Montréal, Montréal, Canada}

\affil[2]{GERAD and School of Mathematics, University of Edinburgh, Edinburgh, United Kingdom}

\affil[3]{Institut de Recherche d’Hydro-Québec, Hydro-Québec, Varennes, Canada}




\date{ }

\maketitle

\begin{abstract}
    With the increasing effects of climate change, the urgency to step away from fossil fuels is greater than ever before. Electric vehicles (EVs) are one way  to diminish these effects, but their widespread adoption is often limited by the insufficient availability of charging stations. In this work, our goal is to expand the infrastructure of EV charging stations, in order to provide a better quality of service in terms of user satisfaction (and availability of charging stations). Specifically, our focus is directed towards urban areas. We first propose a model for the assignment of EV charging demand to stations, framing it as a maximum flow problem. This model is the basis for the evaluation of user satisfaction with a given charging infrastructure. Secondly, we incorporate the maximum flow model into a mixed-integer linear program, where decisions on the opening of new stations and on the expansion of their capacity through additional outlets is accounted for. We showcase our methodology for the city of Montreal, demonstrating the scalability of our approach to handle real-world scenarios. We conclude that considering both spacial and temporal variations in charging demand is meaningful when solving realistic instances.
\end{abstract}





\section{Introduction}\label{sec:introduction}

   Transportation accounts for 28\% of greenhouse gas (GHG) emission in the US~\cite{USEPA} and similarly in the UK (27\%) and in Canada (28\%)~\cite{EEA,GoC2022}. For countries where a large percentage of electricity is generated from renewable sources, as is the case in Canada, studies show that electric vehicles (EVs) are a good alternative to fuel-based ones as a measure to curtail GHG emissions~\cite{Axsen2015,Woo2017}. To boost EV adoption, expansion and improvements to the already existing charging infrastructures must be made. This is because the willingness of car users to opt for an EV is closely linked to the EVs' travel range and the availability of charging stations~\cite{Pevec2020}. The addition of new charging stations can alleviate range anxiety, especially for prospective EV owners~\cite{Carley2013}. As such, Hydro-Qu\'ebec, a publicly owned company responsible for most of the electric grid in the province of Quebec, is investing into more and faster charging stations. In fact, from 2017 to 2022, 1,800 new charging stations were added in the province of Quebec.  Simultaneously, there was a surge in EV purchases, escalating from 3,347 in 2017 to 34,082 in 2022~\cite{StatCan}. This increase in EV purchases is mostly likely influenced by government policies (e.g., \cite{GoC2023,Quebec}), yet it may also be due to the introduction, in urban areas, of  charging stations near homes, workplaces and public areas, as this is known to be a crucial incentive for EV adoption~\cite{Hardman2018}. Homes can sometimes have access to privately owned chargers~\cite{Bailey2015}, but this does not apply to every EV owner. To top it off, public charging infrastructure has been shown to improve EV adoption~\cite{Coffman2017}. This unfortunately leads to the "chicken and egg" dilemma~\cite{Anjos2020}, where investors are only willing to supply more infrastructure if adoption is high, while EV purchases are dependent on widespread charging availability. As such, initial investment must come from governments and public institutions.

The motivation behind this work is to assist decision-makers who need to identify where infrastructure improvements are needed. Concretely, given a set of candidate locations for opening new stations and existing ones, infrastructure owners must decide where to install new stations and determine the number of outlets to be added to new and existent stations. This decision must be driven by the EV users demand. In a urban planning context, it is expected that users mostly carry out intracity trips between home, workplace and public areas, which can serve to identify the areas of infrastructure enhancement. We assume that it is possible to reduce these users' range anxiety by providing them with access to charging stations near these origin and destination places. 

\paragraph*{\textbf{Contributions}} Motivated by the context presented above, this paper tackles the challenge of optimally locating and sizing EV charging stations in urban areas to maximize the satisfied charging demand. Here, charging demand refers to requests from EV users to charge their vehicle at a public station, while satisfying that demand implies the availability in both time and space of a charging station with sufficient capacity for that user. Maximizing the satisfied charging demand is an important tactical planning problem faced by EV infrastructure providers like Hydro-Qu\'ebec, which regularly take decisions on the expansion of their infrastructure to meet the growing charging demand based on current usage. Our first contribution is the formulation of a linear programming model to efficiently evaluate the satisfied demand for a group of existing charging stations. Even though the satisfied demand can be determined from usage data of stations, this model serves two purposes: \emph{(i)} it also allows us to compute the unsatisfied demand and \emph{(ii)} it enables us to evaluate the satisfied demand for any set of stations. Expanding upon this, our second contribution is the integration of location and sizing decisions in the formulation, given a list of candidate locations for new stations, resulting in a mixed-integer linear program which can be used to maximize the satisfied demand. Lastly, our third contribution involves detailing a case study of the island of Montreal. We base our research on real charging session data and origin-destination (OD) trips across Montreal boroughs. With it, we validate the effectiveness of our approach to solve large-scale instances and we conduct an analysis of the solutions it produces.

Our methodology differs from most papers in the literature in three key ways, underlying the novelty of our contributions. 
\begin{enumerate}
    \item We solve the problem of determining the charging demand, i.e., the assignment of the EV users (demand) to stations, by formulating it as a maximum flow problem. Maximum flow problems have the advantage of being solvable efficiently. Importantly, our maximum flow problem is based on~\cite{Ford1958} which is different from the flow-based model commonly found in~\cite{Kuby2005} and other papers about the station location problem. To the best of our knowledge, this is the first maximum flow model of its kind used within the context of the EV station location and sizing problem.
    \item Existing EV station placement methods can handle capacities (e.g., \cite{Upchurch2009}), multiple periods (e.g., \cite{Zhang2017}), already existing infrastructure (e.g., \cite{Yang2018}), large instances (e.g., \cite{Shahraki2015}) or exact solutions (e.g., \cite{Cavadas2015}). While the existing literature usually tackles only one or two of these aspects at a time, this paper integrates all of them in the maximum flow formulation.
    \item Thanks to our partnership with Hydro-Qu\'ebec, we have access to real-world data, including the existing station locations and charging sessions with timestamps and energy consumption for every user. We use this information to generate realistic instances for testing our methodology, and demonstrate our ability to solve large-scale instances with hundreds of stations and the aggregated power demand of thousands of users. 
\end{enumerate}
    
\paragraph*{\textbf{Paper organization}} The paper is organized in the following way. In Section~\ref{sec:literature_paper}, we provide an overview of the existing literature on EV charging infrastructure planning, focusing particularly on station placement. In Section~\ref{sec:model_paper}, we present the linear model for charging station network evaluation in terms of satisfied demand and the mixed-integer program, including station location and sizing decisions. In Section~\ref{sec:experiments}, we describe our case study for the island of Montreal and test our models on realistic instances. Section~\ref{sec:conclusion} concludes the paper and proposes potential future research.

\section{Related Literature}\label{sec:literature_paper}

In this section, we begin with a brief review of the literature pertaining to the optimization of charging infrastructure utilization. Then, we delve into the literature's approaches to estimate charging demand, a crucial element for the optimal placement of charging stations. Subsequently, we discuss different location models, objective functions, intracity and intercity case studies, and temporal modeling considerations. Lastly, we position our work within the reviewed literature.

Research has been conducted on optimizing the existing charging infrastructure, particularly, through charging price decisions aimed at managing the distribution of demand (e.g., \cite{Flath2014}, \cite{Hu2016}, \cite{Moghaddam2019}). However, in our case study of the island of Montreal, prices cannot be changed and the power grid is prepared to handle even the most severe winter day. Therefore, we focus our review on charging station placement.

Decisions on the expansion and opening of charging stations requires the knowledge of its potential use. Thus, the estimation of charging demand is important, as it indicates when and where the demand for charging originates. The most common method is based on the use of OD data to subsequently model how users move between locations. This kind of data often comes from surveys (e.g., \cite{Baouche2014}, \cite{Zhang2015}, \cite{Cavadas2015}). This is the approach adopted in this paper. Nonetheless, it is important to note that our data only covers travels between boroughs (i.e., it is not granular) and it is not exclusive to EVs. Hence, we complement the demand estimation with other available data.

The location modelling usually falls into one of two categories: node-based or flow-based~\cite{Upchurch2010}. In the node-based approach, either a list of candidate locations is provided and the goal is to maximize the coverage (e.g., \cite{Frade2011}, \cite{Tu2016}, \cite{Yang2018}), or population nodes are used as candidate locations and the goal is to satisfy all the demand at minimum cost, i.e., cost of opening stations~\cite{Zhang2015,Li2016,Xie2018, Bouguerra2019}. For the flow-based modelling, flow is assigned OD pairs, and facilities (charging stations in this context) must capture as much flow as possible. This is another variant of maximum coverage proposed by~\cite{Hodgson1990}. Using the flow-based modelling, Kuby and Lim~\cite{Kuby2005} are the first to propose the Fuel Refuelling Location Problem (FRLP) which seeks to locate a fixed number of refuelling stations on a network so as to maximize the total flow volume refuelled.  In our work, since we consider intracity travels, and hence, short trips, we do not consider the routing of EVs. We define a maximum flow problem in the sense of \cite{Ford1958} for the location modelling. The key difference with the FRLP is that we treat flow as a variable rather than a parameter. In the FRLP each OD is assigned a flow volume on the shortest path between the origin and destination. A binary variable is then used to identify whether each flow volume is present or not when maximizing the objective. This is fundamentally different from our approach, since we view flow as a variable which can enter and leave both OD pairs and stations using flow constraints.

In the literature, the objective of the charging station placement problems varies significantly, but it is often closely related to the location modelling choice. Flow-based models often maximize the total amount of flow in the network (e.g., \cite{Kuby2005}, \cite{Capar2013}, \cite{Chung2015}, \cite{Kadri2020}),  which is also our case. For node-based modelling,  different objectives have been used, such as maximization of the satisfied EV demand~\cite{Cavadas2015,Tu2016,Yang2018} and minimization of the costs~\cite{Zhang2015,Li2016,Yang2017,Xie2018,Bouguerra2019,Zhong2022,Filippi2023}.

A key aspect of the EV charging station placement problem is whether it  concerns the intercity or intracity context. The intercity case focuses on long distance travels between cities, with users potentially charging once or more during a trip (e.g., \cite{Chung2015}, \cite{Li2016}, \cite{Xie2018}). The intracity case focuses on a city, with users typically charging near homes, workplaces or public areas~\cite{Hardman2018}. The works in \cite{Frade2011,Baouche2014,Cavadas2015} are all examples of works on intracity problems. To the best of our knowledge, Anjos, Gendron and Joyce-Moniz~\cite{Anjos2020} are the only authors to handle both intra and intercity cases simultaneously. 

Some works include a time component to the EV charging station placement problem.  Notably, time periods have been modelled with two different goals: to consider strategical planning, where periods can represent years, and to consider tactical planing, where periods can represent hours. If time is modelled over a certain number of years, the goal is to focus on the adoption of EVs or the evolution of the EV infrastructure~\cite{Chung2015,Li2016,Zhang2017, Xie2018,Anjos2020,Lamontagne2023}. If time is modelled over a range of hours, the goal is to reflect high and low demand over certain times and evaluate the infrastructure service quality~\cite{Frade2011,Cavadas2015,Tu2016,Filippi2023}.

In this paper, we propose a flow-based (in the sense of \cite{Ford1958}) mixed-integer linear program (MILP) for the EV charging station placement and sizing problem. Baouche \emph{et al.}~\cite{Baouche2014} investigated a case study of the city of Lyon for their intracity model. They also used OD data to estimate demand, yet the approach is fundamentally different from ours, with their optimization model ensuring that all demand is covered at minimum cost, and without complementing the OD data with EV session data. Filipi \emph{et al.}~\cite{Filippi2023} emphasized the importance of accounting for spatial and temporal variations in demand, which aligns with the considerations in our study. They adopted a node-based demand model and focused on minimizing installation costs and customers travel distance, subject to satisfying all the demand (which can be assigned to any opened station). This contrasts with our approach in two key ways: firstly, we focus on maximizing the satisfied energy demand; secondly, we use OD data and we limit the feasibility of the charging stations to points near the origin or destination, rather than node-based demand modelling. Lamontagne \emph{et al.}~\cite{Lamontagne2023} proposed a MILP for the maximization of EV adoption in the long-term. Their model does not consider the stations' capacity, which is a crucial factor for our tactical problem, maximizing satisfied demand. Cavadas, H. A. Correia and Gouveia~\cite{Cavadas2015} considered an intracity case study as well as a temporal dimension along with station capacities. Their work differentiates from ours as they use a node-based model rather than flow-based, aiming to minimize walking distance, and using a predefined number of outlets per station. Finally, the direct solving of the various mixed-integer linear programs proposed for the EV placement problem has encountered the issue of scalability (e.g., \cite{Zhang2017,Anjos2020,Lamontagne2023}). However, by leveraging on the maximum flow model for the estimation of the satisfied demand, we are able to solve our MILP for instances based on the real Montreal demand and existing EV infrastructure.

\section{Mathematical Formulation}\label{sec:model_paper}
\subsection{Problem Statement}

Our problem involves determining the optimal location and sizing (number of outlets) of EV charging stations in an urban context. The urban area under study can possess existing stations, but it is not required for the correctness of our model. The decision-maker's goal is to maximize the satisfied daily (charging) demand subject to a budget constraint for the infrastructure expansion costs. Certainly, to address this problem, it is crucial to model how current EV users utilize the available charging stations, either existing or newly installed. Hence, we next describe the available information about EV users and the assumptions made in our work.

Since we consider the urban case, we expect EV users to travel between home and work, home and childcare, home and leisure areas, and so on, which are relatively short distances within the range of EVs. Therefore, we can assume that they do not charge along a path but rather at its origin or destination. Hence, the problem of determining how the charging demand is spread over the available stations becomes a matching problem, where we aim to match EV users to stations close to their origin or destination.  Maximizing the number of matchings is equivalent to determining the maximum demand that can be satisfied. Given that, in our case study, users have access to an app providing in real-time  the information about station occupancy (\href{https://lecircuitelectrique.com/en/mobile-app}{\emph{Le Circuit électrique}}\footnote{\url{https://lecircuitelectrique.com/en/mobile-app/}}), it is reasonable to optimize the assignment with this objective function.

Another important aspect of our problem is the consideration of time. Over a day, EV users do not necessarily travel and charge at the same time, nor do charging sessions have the same duration. For instance, we should expect peaks of demand in the evenings in residential areas, and significant charging duration differences between level 2 and level 3 charging stations \cite{Morrissey2016}. Therefore, we discretize the day into a finite number of periods over which the demand varies, and we consider the assignment of users to stations for each of these periods.

In our case study, we have access to the origin-destination matrix for the urban area under investigation, along with charging session data for existing stations.
\subsection{Linear Model: Assigning Demand to Stations}

In this section, we describe our framework to determine the assignment of EV charging demand to stations. To this end, we first provide a graph modeling and then, a linear programming formulation. Our notation is summarized in Table~\ref{table:notation}; the elements corresponding to new stations and outlets are only used in the next section. 

\begin{table}[h!] \footnotesize 
\caption {Notation}
\label{table:notation}
\begin{tabular}{ccl} 
\hline
\textbf{Type}                  & \textbf{Notation} & \multicolumn{1}{c}{\textbf{Description}}                    \\ 
\hline
\multirow{9}{*}{\textbf{Sets}} & $T$        & Set of time periods \\
                               & $V$        & Set of vertices $\{1,2, ..., N\}$ where 1 is the source and \textit{N} is the sink \\
                               & $O$        & Subset of vertices representing OD pairs   \\ 
                               & $S_1$      & Subset of vertices representing existing stations  \\ 
                               & $S_2$      & Subset of vertices representing candidate locations\\ 
                               & $L$        & Subset of arcs representing the charging flow demand of users per OD pairs \\ 
                               & $M$        & Subset of arcs between OD pairs and stations \\ 
                               & $R_1$      & Subset of arcs representing the charging flow supply at an existing station \\ 
                               & $R_2$      & Subset of arcs representing the charging flow supply at a candidate location \\ 
\hline
\multirow{10}{*}{\textbf{Parameters}} & $A^t_e$    & Charging  flow demand of an OD in period $t$ for $e \in L$\\ 
                                     & $C_e$      & Charging flow supply at an existing station $e \in R_1$  \\ 
                                     & $I^{(2)}_e, I^{(3)}_e$      & Cost of installing an outlet to a new level 2 or 3 station in location $e \in R_2$ \\
                                     & $J^{(2)}_e, J^{(3)}_e$      & Cost of building a new level 2 or 3 station in location $e \in R_2$ \\
                                     & $K_e$              &   Cost of adding an outlet to an existing station $e \in R_1$  \\  
                                     & $G$                & Budget   \\  
                                     & $P_e$              & Amount of charging flow supply for a single outlet at an existing station $e \in R_1$ \\  
                                     & $Q^{(2)}, Q^{(3)}$        & Amount of charging flow supply for a single outlet at a new level 2 or 3 station   \\   
                                     & $Y^{(2)}, Y^{(3)}$        & Maximum number of outlets in a level 2 or 3 station \\  
                                     & $Y_e$                 & Maximum number of outlets in location $e \in R_1$ \\
\hline
\multirow{7}{*}{\textbf{Variables}} & $a^t_e$       & Amount of charging flow demand generated by an OD in period $t$ for edge $e \in L$ \\
                                    & $b^t_e$       & Amount of charging flow going from an OD to a station in period $t$ for edge $e \in M$  \\ 
                                    & $c^t_e$       & Amount of charging flow supply going through an existing station in period $t$ for edge $e \in R_1$  \\ 
                                    & $d^t_e$             & Amount of charging flow supply going through a candidate location in period $t$ for edge $e \in R_2$ \\
                                    & $x_e$               & Number of outlets to add to an existing station $e \in R_1$ \\
                                    & $y^{(2)}_e, y^{(3)}_e$       & Number of outlets of a new level 2 or 3 station $e \in R_2$   \\ 
                                    & $z^{(2)}_e, z^{(3)}_e$       & Binary variable indication whether or not to build a new level 2 or 3 station  $e \in R_2$   \\ 
\hline
\end{tabular}
\end{table}

\paragraph*{\textbf{Graph transformation}}
Let us go into more detail in the description of our assignment problem as it consists of a crucial building block for our methodology. We use a bipartite graph to describe potential matches (assignments). The left side of the bipartite graph is composed of EV users and the right side of stations. Instead of using individual EV users as vertices, it is more efficient to group them based on their trips: users with the same origin-destination (OD) pair are grouped together. The edges must represent feasible stations for each OD pair. To decide if a station is feasible for a given OD pair, we define a parameter $R$ which describes the maximum radius around the origin or the destination of an OD pair. If a station is within either radius, then it is feasible for that OD pair and an edge is created.

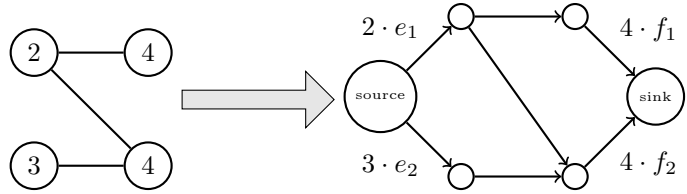
\begin{figure}[!h]
\centering

\begin{tikzpicture}[node distance={15mm}, thick, main/.style = {draw, circle}] 
\node[main] (1) {$2$}; 
\node[main] (2) [below of=1] {$3$}; 
\node[main] (3) [right of=1] {$4$};
\node[main] (4) [right of=2] {$4$};

\draw (1) -- (3);
\draw (1) -- (4);
\draw (2) -- (4);
\end{tikzpicture}
\begin{tikzpicture}[baseline=-22mm]
    \node[single arrow,draw=black,fill=black!10,minimum height=2cm,shape border rotate=0] at (0,-1) {};
\end{tikzpicture}
\begin{tikzpicture}[node distance={15mm}, thick, main/.style = {draw, circle}] 
\node[main] (5) [below right of=3] {\tiny source}; 

\node[main] (6) [above right of=5] {};
\node[main] (7) [below right of=5] {}; 
\node[main] (8) [right of=6] {};
\node[main] (9) [right of=7] {};
\node[main] (10) [below right of=8] {\tiny sink };

\draw[->] (5) -- node[midway, above left] {$2 \cdot e_1$} (6);
\draw[->] (5) -- node[midway, below left] {$3  \cdot e_2$} (7);
\draw[->] (6) -- (8);
\draw[->] (6) -- (9);
\draw[->] (7) -- (9);
\draw[->] (8) -- node[midway, above right] {$4  \cdot f_1$} (10);
\draw[->] (9) -- node[midway, below right] {$4  \cdot f_2$} (10);

\end{tikzpicture} 
\caption {Converting a bipartite graph to a flow graph}
\medskip
\footnotesize
Example: The bipartite graph on the left has two ODs, each with 2 and 3 EV users, and two stations, each with 4 outlets; the edges represent the feasible stations.  On the right, we have the transformed flow graph, where in some of the arcs we have their maximum flow capacity related to EV demand and station supply; the conversion factors $e_1$ and $e_2$ map users to flow, while $f_1$ and $f_2$ map outlet supply to flow.
\label{fig:flow_Graph}
\end{figure}

Issues with this maximum matching approach are that there can exist more than one user per OD pair and a station can charge more than one user at a time; note that a station can have more than one outlet. To fix this, we convert the bipartite graph into a flow graph, and the matching problem into a maximum flow problem. To do so, the edges of the bipartite graph are transformed into arcs from the vertices in set $O$ to the vertices in set $R_1$, a source vertex and a sink vertex are introduced, an arc from the source to each element in set $O$ is added, and an arc from each element in set $R_1$ to the sink is added. Finally, to define a maximum flow problem over the resulting graph, restrictions regarding the amount of flow that can pass through each arc and the cost of using those arcs must also be defined. In our problem, there is no cost associated with the use of the arcs. It would be possible to add a cost based on the distance between an origin or destination and a station but this is out of the scope of this paper, since we assume that stations are within a short walking distance. On the other hand, we do define maximum flow capacities to the arcs between the source and the elements in set $O$, and between elements in set $R_1$ and the sink. For the former, the maximum flow capacity represents, at a given period, the amount of charging flow demand for every user travelling on each OD. For the latter, the maximum flow capacity on the arc from a station to the sink represents the maximum amount of charging flow supply available at that station within a period. See Figure~\ref{fig:flow_Graph} for an illustration. A key aspect here is that the charging flow is relative to a period. As such, the (flow) graph can be replicated for a given set of periods $T$, where the graph remains the same but the arcs' maximum capacities can change between periods. Specifically, the maximum flow capacities on the arcs from the source to the OD pairs may vary, allowing for the representation of fluctuating number of users travelling on OD pairs at different times of the day. In this way, if we determine the maximum flow from the source to the sink of our graph over a finite time horizon (in our case study, 24 hours), we determine the maximum (daily) charging demand that can be satisfied by the current infrastructure.

Note that, in order to identify if a potential station location is interesting, we can add a list of candidate locations to the flow graph, following a similar approach as with existing stations. For instance, when generating the flow graph, it is possible to encounter OD pairs with no arc connected to a station. The demand from such an OD is referred to as \emph{impossible demand}. Therefore, it would make sense to have in the list of potential new stations one or more locations close to the said OD pair; if we open at least one of these stations, it would guarantee an increase in the overall satisfied demand.

\paragraph*{\textbf{Formulation}} We are ready to  provide the linear program corresponding to the maximum flow of the described graph. 

\begin{figure}[!h]
\centering
\begin{tikzpicture}[node distance={15mm}, thick, main/.style = {draw, circle}] 
\node[main, blue] (5) [below right of=3] {\tiny source}; 

\node[main, orange] (6) [above right of=5] {};
\node[main, orange] (7) [below right of=5] {}; 
\node[main, green] (8) [right of=6] {};
\node[main, green] (9) [right of=7] {};
\node[main, pink] (10) [below right of=8] {\tiny sink};

\draw[->, red] (5) -- (6);
\draw[->, red] (5) -- (7);
\draw[->, violet] (6) -- (8);
\draw[->, violet] (6) -- (9);
\draw[->, violet] (7) -- (9);
\draw[->, brown] (8) -- (10);
\draw[->, brown] (9) -- (10);

\end{tikzpicture} 
\caption{Flow model notation}
\medskip
\footnotesize
From left to right, \textcolor{blue}{\rule{1ex}{1ex}}~source (vertex $1$), \textcolor{red}{\rule{1ex}{1ex}}~set $L$, \textcolor{orange}{\rule{1ex}{1ex}}~set $O$, \textcolor{violet}{\rule{1ex}{1ex}}~set $M$, \textcolor{green}{\rule{1ex}{1ex}}~set $S_1 \cup S_2$, \textcolor{brown}{\rule{1ex}{1ex}}~set $R_1 \cup R_2$, \textcolor{pink}{\rule{1ex}{1ex}}~sink~(vertex $N$).
\label{fig:notationflow}
\end{figure}
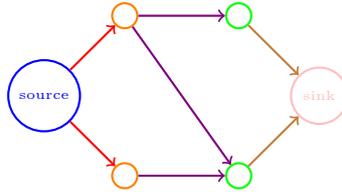

Figure~\ref{fig:notationflow} provides a visual summary of the notation used for the sets of arcs (continuation of the example of Figure~\ref{fig:flow_Graph}).

Our maximum flow problem is the following linear program:
\begin{subequations}
\label{problemlinear}
    \begin{align}
    \max_{a,b,c} & \sum_{t \in T}\sum_{e \in L} a_e^t \label{lin:objective}\\
    s.t. \  & a_{(1,v)}^t = \sum_{e \in M:e=(v,i)} b_e^t &\forall v \in O,\forall t \in T \label{lin:flow1}\\
    & \sum_{e \in M: e=(i,v)} b_e^t = c_{(v,N)}^t &\forall v \in S_1,\forall t \in T \label{lin:flow2}\\
    & 0 \le a_e^t \le A_e^t & \forall e \in L,\forall t \in T \label{lin:constraint_a}\\
    & 0 \le b_e^t & \forall e \in M,\forall t \in T \label{lin:constraint_b}\\
    & 0 \le c_e^t \le C_e & \forall e \in R_1,\forall t \in T. \label{lin:constraint_c}
    \end{align}
\label{Program:linear}
\end{subequations}
The objective function \eqref{lin:objective} is the sum of charging flow leaving the source. Since the flow Constraints~\eqref{lin:flow1} and~\eqref{lin:flow2}  guarantee that the amount of flow reaching the sink is equal to the amount leaving the source, this objective function is equivalent to the total amount of charging flow in the graph.  Indeed, the flow Constraints \eqref{lin:flow1} and \eqref{lin:flow2} ensure that the amount of charging flow entering into the OD pairs is the same amount leaving and the amount of charging flow entering the stations is the same amount leaving, respectively. Constraints~\eqref{lin:constraint_a} limit the flow from the source to each OD pair to the charging demand for that specific OD pair within a given period. This charging demand is equal to the number of EV users travelling on that OD pair multiplied by the charging flow demand per user. For example, let a flow $f$ be the average kilowatt per second consumption of an EV and $p$ the length of a period in seconds, the charging flow demand of a user is $f \cdot p$. If we multiply that value by the number of EV users in an OD, we obtain the complete requested charging flow for that OD. Constraints~\eqref{lin:constraint_c} limit the amount of charging flow supply that can be provided at each station. 

\subsection{Mixed-Integer Model: Placing Stations and Outlets}
In the previous section, we described our model for the assignment of the demand to stations. Next, our objective is to introduce new stations  and outlets that remain available throughout all periods, with the aim of diminishing the existing unsatisfied demand (or, equivalently, increase the satisfied demand). It is worth noting that our approach assumes that all stations and outlets are built from the beginning, rather than gradually over time. We integrate into Program~\eqref{problemlinear} the decisions related to the opening of new stations and addition of outlets, leading to the following mixed-integer program:
\begin{subequations}
    \label{problemmixed-integer}
    \begin{align}
    \max_{a,b,c,d,x,y,z} & \sum_{t \in T}\sum_{e \in L} a_e^t \\
    s.t. & \sum_{e \in R_2} \left(I_e^{(2)} y_e^{(2)} + J_e^{(2)} z_e^{(2)} + I_e^{(3)}  y_e^{(3)} + J_e^{(3)} z_e^{(3)}\right)+ \sum_{e \in R_1} K_e x_e\le G \label{mip:budget}\\
    & \eqref{lin:flow1}-\eqref{lin:constraint_b} \nonumber\\
     & \sum_{e \in M: e=(i,v)} b_e^t = d_{(v,N)}^t \hspace*{5cm}  \forall v \in S_2,\forall t \in T \label{mip:flow2} \\
    & 0 \le c_e^t \le C_e + P_e x_e  \hspace*{5cm}   \forall e \in R_1,\forall t \in T \label{mip:constraint_c}\\
    & 0 \le d_e^t \le Q^{(2)} y_e^{(2)} + Q^{(3)}  y_e^{(3)} \hspace*{3.9cm}  \forall e \in R_2,\forall t \in T \label{mip:constraint_d}\\
    & x_e \le Y_e - \frac{C_e}{P_e}  \hspace*{6.9cm} \forall e \in R_1 \label{mip:outlet_x}\\
    & y_e^{(2)} \le Y^{(2)} z_e^{(2)} \hspace*{6.9cm} \forall e \in R_2 \label{mip:outlet_y}\\
    & y_e^{(3)} \le Y^{(3)} z_e^{(3)}  \hspace*{6.9cm} \forall e \in R_2 \label{mip:outlet_y2}\\
    & z_e^{(2)} + z_e^{(3)} \le 1  \hspace*{6.9cm} \forall e \in R_2 \label{mip:constraint_z}\\
    & x_e \in \mathbb{N}  \hspace*{7.8cm} \forall e \in R_1 \label{mip:constraint_x}\\
    & y_e^{(2)}, y_e^{(3)} \in \mathbb{N}, z_e^{(2)}, z_e^{(3)} \in \{0, 1\}  \hspace*{4.6cm} \forall e \in R_2. \label{mip:constraint_yz}
    \end{align}
\label{Program:mip}
\end{subequations}
The objective function is the same as before. The constraint \eqref{mip:budget}  enforces  the costs of new stations and outlets to be below a given budget. Constraints~\eqref{mip:flow2} are the same as Constraints~\eqref{lin:flow2} but for candidate stations instead. Constraints~\eqref{mip:constraint_c} adapt Constraints~\eqref{lin:constraint_c} to account for newly added outlets. The number of new outlets is multiplied by a factor $P_e$ to convert them into flow. Constraints~\eqref{mip:constraint_d} limit the maximum amount of flow that can travel to candidate stations by the amount of new level 2 or level 3 outlets. Constraints~\eqref{mip:outlet_x} limit the maximum amount of newly added outlets to an existing station by subtracting the already existing number of outlets from the maximum number possible. The parameters $Y_e$ and the variables $x_e$ are relative to the level of the station $e$ they are associated with. If a station is level 2, only level 2 outlets can be installed and the same applies to level 3. Constraints~\eqref{mip:outlet_y} and \eqref{mip:outlet_y2} limit the number of outlets according to whether a new level 2 or level 3 station is built. Due to Constraints~\eqref{mip:constraint_z}, a new station can only be level 2 or level 3 but not both. In the rare case where there are level 2 and level 3 outlets at an existing station, that station is considered as two separate stations in the same location. Constraints~\eqref{mip:constraint_x} and \eqref{mip:constraint_yz} set the domains for variables $x, y, z$.

\section{Computational Experiments}\label{sec:experiments}

In this section, we aim to validate the use of our linear program to estimate station demand and the efficiency of solving our mixed-integer model for real-world instances. Hence, we start in Section~\ref{subsec:MontrealCase} by detailing our case study of the island of Montreal. Then, in Section~\ref{subsec:AssignUsers}, we show experimental results of our linear model when it comes to matching users to existing stations. Finally, in Section~\ref{subsec:AddingStations}, we provide experimental results for solving our mixed-integer program, modelling the addition of new stations and outlets. All experiments are run on an Intel i7-10700F CPU @ 2.90 GHz with 8 cores and 16 GB of RAM. We use CPLEX Optimization Studio V22.1.0 on a single thread per instance and a 30 minute time limit.

\subsection{Montreal Case}\label{subsec:MontrealCase}

\begin{figure}  
\centering  

\subfigure[Level 2]  
{  
\includegraphics[scale=0.45]{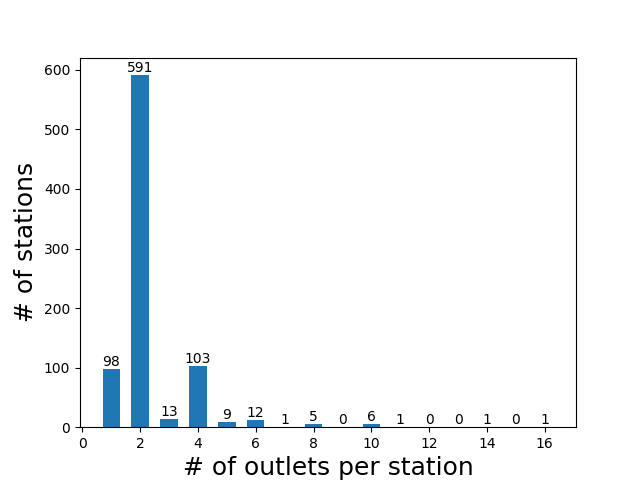}
}  
\subfigure[Level 3]  
{  
\includegraphics[scale=0.45]{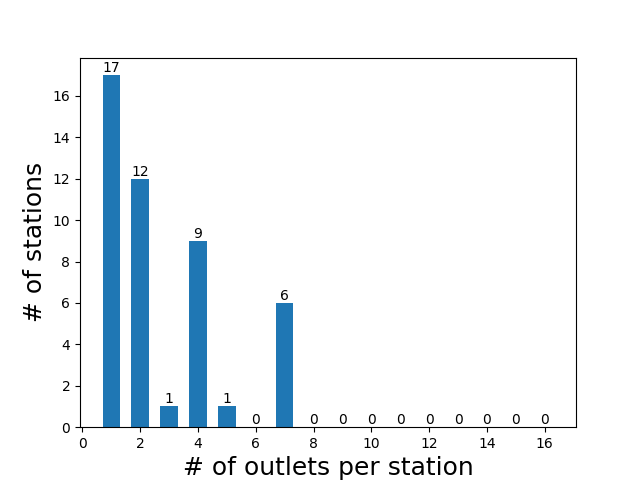}
}

\caption{Distribution of the number of outlets per station}
\label{fig:outlethistogram}
\end{figure}

\paragraph*{\textbf{Data}}
We focus our experiments on the island of Montreal. For this case study, we obtained data about the location, level and number of outlets for existing public stations of the \emph{Le Circuit électrique} and the time, duration and average kilowatts per second (kW/s) for charging sessions at these stations. In our data, there are 841 level 2 stations and 41 level 3 stations on the island. The maximum number of outlets within a station is 16 and 7 for level 2 and level 3 stations, respectively. Figure \ref{fig:outlethistogram} provides the distribution of outlets per station. 

    
    

\begin{figure}  
\centering  

\subfigure[Level 2]  
{  
\begin{tikzpicture}[scale=0.8]
  \begin{axis}
    [
    ytick={1,2,3,4},
    yticklabels={0h-6h, 6h-12h, 12h-18h, 18h-24h},
    ]
    \addplot+[
    boxplot prepared={
      median=36.986647099992576,
      upper quartile=41.73775323829673,
      lower quartile=29.84551632287221,
      upper whisker=79.80767288228044,
      lower whisker=0.0
    },
    ] coordinates {};
    \addplot+[
    boxplot prepared={
      median=21.214313378505245,
      upper quartile=24.499076552719407,
      lower quartile=18.932434529412376,
      upper whisker=61.29550424803081,
      lower whisker=0.0
    },
    ] coordinates {};
    \addplot+[
    boxplot prepared={
      median=19.934688865399675,
      upper quartile=23.18482216359185 ,
      lower quartile=17.722731908076514,
      upper whisker=84.24986652429258,
      lower whisker=0.0
    },
    ] coordinates {};
    \addplot+[
    boxplot prepared={
      median=21.643616328026835,
      upper quartile=24.051223856310353,
      lower quartile=19.820234930397124,
      upper whisker=100.0,
      lower whisker=0.0
    },
    ] coordinates {};
  \end{axis}
\end{tikzpicture}

}  
\subfigure[Level 3]
{  

\begin{tikzpicture}[scale=0.8]
  \begin{axis}
    [
    ytick={1,2,3,4},
    yticklabels={0h-6h, 6h-12h, 12h-18h, 18h-24h},
    ]
    \addplot+[
    boxplot prepared={
      median=23.525968066988498,
      upper quartile=25.51334755841905,
      lower quartile=21.744155055699316,
      upper whisker=33.262649460696153,
      lower whisker=0.0
    },
    ] coordinates {};
    \addplot+[
    boxplot prepared={
      median=26.17619080342869,
      upper quartile=27.640433169462647,
      lower quartile=25.401388421919785,
      upper whisker=42.567893299051907,
      lower whisker=0.0
    },
    ] coordinates {};
    \addplot+[
    boxplot prepared={
      median=25.619836732935924,
      upper quartile=26.836278602851804  ,
      lower quartile=24.355330844302003,
      upper whisker=100.0,
      lower whisker=15.702982349360928
    },
    ] coordinates {};
    \addplot+[
    boxplot prepared={
      median=24.122496449115993,
      upper quartile=25.919283906305535,
      lower quartile=22.8616210009872,
      upper whisker=39.30580106058172,
      lower whisker=0.0
    },
    ] coordinates {};
  \end{axis}
\end{tikzpicture}

}

\caption{Distribution of the daily average percentage of energy supplied per station per 6-hour period}
\label{fig:demandperperiod}
\end{figure}
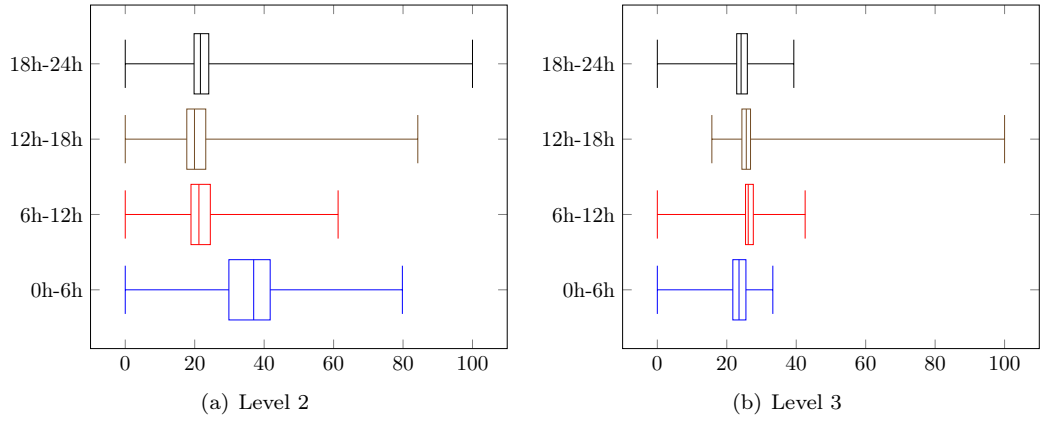

 We tested two sets of periods: \emph{(1)} a single time horizon of 24 hours and \emph{(2)} the same horizon discretized into 6-hour periods. The goal is to find the impact of relaxing the demand over 24 hours, i.e., of assuming that the demand can be satisfied at any moment of the day. Figure \ref{fig:demandperperiod} shows that for level 2 stations, a substantial percentage of the energy provided is to users who leave their vehicles to charge overnight. However, level 3 stations are used uniformly throughout the day as they are more expensive and thus less popular for overnight charging. As such, accounting for fluctuations of the charging demand over the day will result in more accurate modeling of the satisfied demand, which is expected to better inform the placement and sizing of stations. Indeed, using a 24-hour relaxation should overestimate the satisfied demand in comparison with the same discretized horizon, since users can be forced to charge at inconvenient times.

\begin{figure}
    \centering
    \includegraphics[scale=0.5]{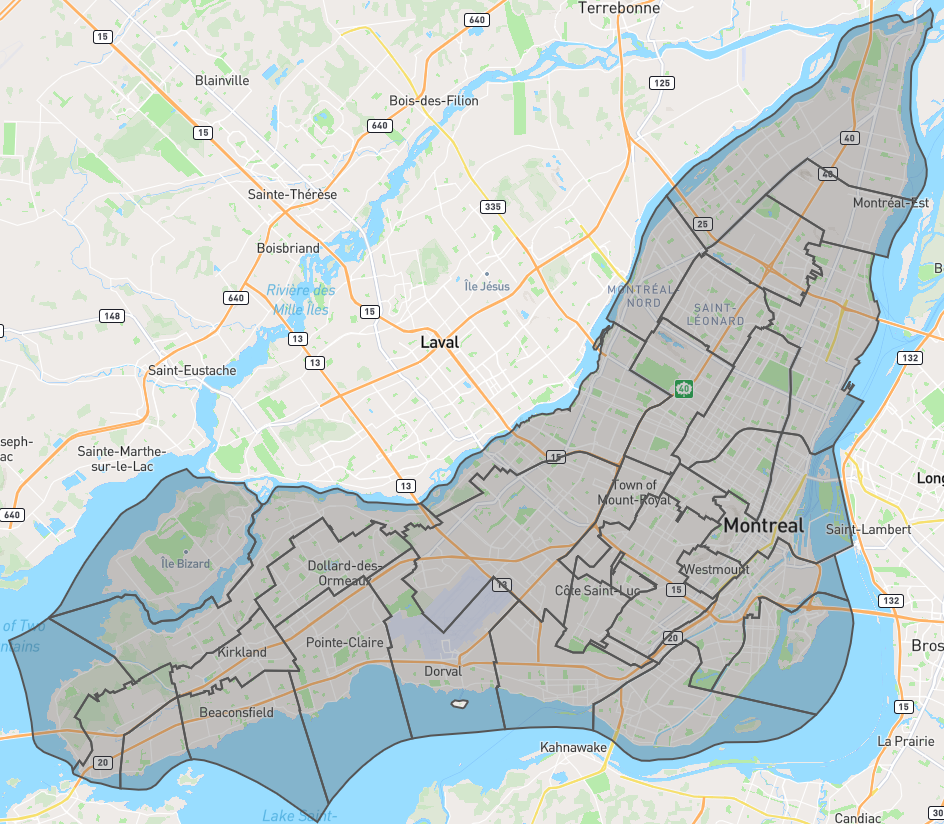}
    \caption{Map of Montreal boroughs}
    \label{fig:montrealmap}
\end{figure}

We also used the publicly available data of the \href{https://www.artm.quebec/planification/enqueteod/}{2018 OD survey}\footnote{\url{https://www.artm.quebec/planification/enqueteod/}} of the Montreal region collected by the ARTM. This data provides the average number of trips between each pair of Montreal boroughs within a day. These trips are not limited to EVs, but can be used as a reasonable approximation of users' movements. We only consider the boroughs within Montreal which leaves us with 32 different boroughs (see Figure \ref{fig:montrealmap}).

\paragraph*{\textbf{Generation of instances}} Based on the described data, we now detail the process used to generate instances, namely, the sets and parameters of Table~\ref{table:notation}. To begin, we need to decide on two parameters: $R$ and $W$. Recall that the parameter $R$ indicates the maximum distance a user is willing to walk between their origin or destination and a charging station. In our framework, each OD pair could have its own radius, however, we use the same radius for all of them since we do not dispose of individual user information related to their willingness to walk. We limit our radius between 400 and 700 meters based on research done on the acceptable walking distance for public transit stops and stores~\cite{Yang2012,Millward2013,Gunn2017,Sugiyama2019}. The parameter $W$ stands for the number of randomly generated (latitude, longitude) points in an instance. More concretely, these randomly generated points serve to create OD pairs, where each point is both an origin and a destination. We generate points relative to the density of EV users charging in each borough (session data). In principle, we want to have a number of points $W$ capable of covering the entirety of the urban area relative to the radius $R$. In practice, however, this would require far too many points or an unrealistically large radius. As such, we try a different number of (randomly generated) points to test the sensibility of our model. In the rest of this section, we provide a simple example to explain each step of our instance generation process and conclude with a summary of the combination of parameters used to generate instances.

\begin{figure}
    \centering
    \includegraphics[scale=1.5]{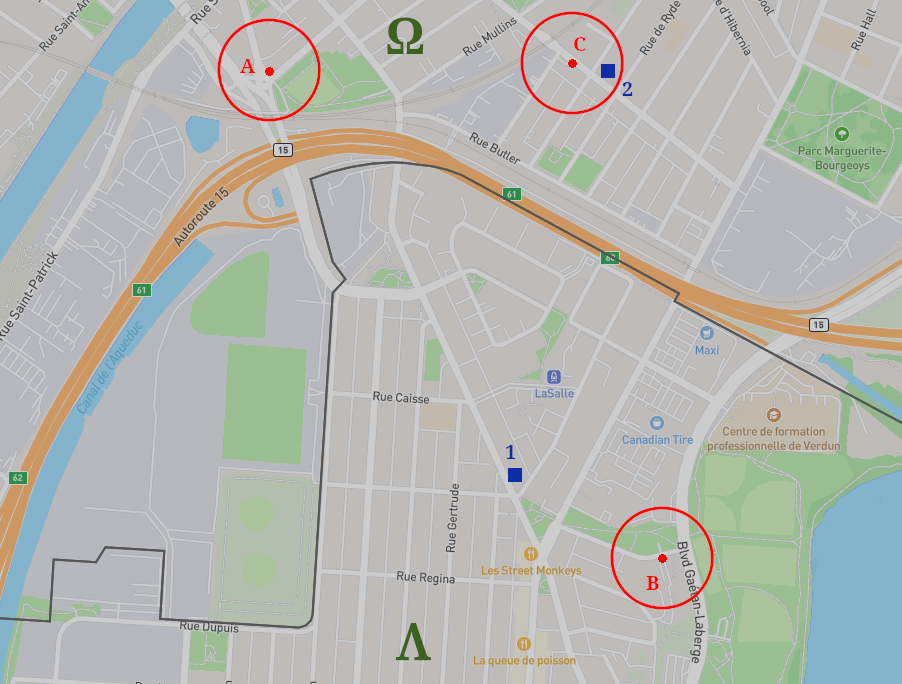}
    \caption{Example of an instance}
    \label{fig:example_map}
\end{figure}

\begin{figure}[!h]
\centering

\begin{tikzpicture}[node distance={15mm}, thick, main/.style = {draw, circle}] 
\node[main] (1) {$AB$}; 
\node[main] (2) [below of=1] {$AC$}; 
\node[main] (3) [below of=2] {$BC$};

\node[main] (7) [right of=1] {$1$};
\node[main] (8) [below of=7] {$2$};

\draw (2) -- (8);
\draw (3) -- (8);

\end{tikzpicture}
\caption {The bipartite graph representation of Figure \ref{fig:example_map}}
\label{fig:example_bipartite}
\end{figure}
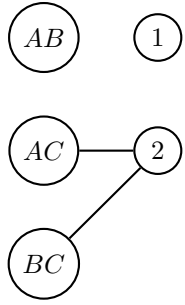

Figure \ref{fig:example_map} is a simple map with two borrows ($\Omega$ and $\Lambda$), containing three randomly generated points in red and two stations represented by the blue squares. Each point has the same radius $R$. We can convert Figure \ref{fig:example_map} into a bipartite graph. Figure \ref{fig:example_bipartite} gives a visual representation of the transformation.

Our model provides the option to use any kind of data as the flow. We use kW since our dataset contains the kW/s per charging session. To calculate the maximum flow capacity of existing stations, we use the session data to estimate the average kW/s per session for each outlet. For each station, we sum up the outlets' kW/s per session to get the maximum kW/s per station. This value can be multiplied by the length of the periods to obtain the maximum flow capacity per period of a station. 

To calculate the $C_e$ parameters of our example, we take the average of kW/s per session from Table \ref{table:example_data} for each station: 5 kW/s for station 1 and 4    kW/s for station 2. If we let our time horizon be a single 24-hour period, then the total flow supply is $5 \cdot 3600 \cdot 24 = 432000$~kW and $4 \cdot 3600 \cdot 24 = 345600$~kW for stations 1 and 2, respectively. For our example, we assume that each station only has one outlet. 

\begin{table}[h!] \footnotesize
\caption {Example of session data}
\centering
\label{table:example_data}
\begin{tabular}{||c c c||} 
 \hline
 Station & Duration (s) & kW/s \\ [0.5ex] 
 \hline\hline
 1 & 110 & 5 \\ 
 \hline
 2 & 100 & 3 \\
 \hline
 2 & 10 & 5 \\
 \hline
\end{tabular}
\end{table}

To calculate the maximum flow capacity per OD (i.e., the demand), we start by calculating the daily average kW that is being used within a period for each station. We sum this supplied kW for each borough; this gives us the amount of supply per borough. From there, we need to convert the supply into its original demand per borough. Here, we combine the use of session data with the Montreal OD data to calculate the percentage of EV users travelling between two boroughs. This forms a set of linear equations: $\displaystyle\sum_{i \in H} q_i p_i^j = r_j$ $\forall j \in H$ where $H$ is the set of Montreal boroughs. In this set of equations, $r_j$ represents the total amount of supply at each station within borough $j$. The coefficient $p_i^j$ is the percentage of people travelling from borough $i$ to $j$. Our variable is $q_i$ which is the total amount of demand in borough $i$. To calculate the amount of demand on each OD, we simply compute $q_i \cdot p_i^j + q_j \cdot p_j^i$ for each borough $i$ and $j$. We assume that most people leave in the morning and come back at night, resulting in bidirectional flow. This assumption is realistic, as most intracity trips occur between home, work and public places. If we want to account for unidirectional trips (which would allow to enforce users to charge only at the origin or only at the destination), we would simply need to duplicate each OD vertex. We distribute the previously computed demand uniformly between each OD pair with the same origin and destination. We also account for trips within the same borough by computing $q_i \cdot p_i^i$.

\begin{table}[h!] \footnotesize
\caption {Example of an OD matrix}
\label{table:od_matrix}
\centering
\begin{tabular}{| c | c | c |} 
 \hline
  & $\Omega$ & $\Lambda$ \\
 \hline
  $\Omega$ & 50\% & 50\% \\ 
 \hline
 $\Lambda$ & 25\% & 75\% \\
 \hline
\end{tabular}
\end{table}

\begin{figure}
    \centering
    \begin{tikzpicture}[node distance={20mm}, thick, main/.style = {draw, circle}] 
    \node[main] (1) {\tiny source}; 
    
    \node[main] (3) [right of=1] {$AC$}; 
    \node[main] (2) [above of=3] {$AB$};
    \node[main] (4) [below of=3] {$BC$}; 
    \node[main] (5) [right of=2] {$1$};
    \node[main] (6) [right of=3] {$2$};
    \node[main] (8) [right of=6] {\tiny sink};

    \draw[->] (1) -- node[midway, above left] {$175$} (2);
    \draw[->] (1) -- node[midway, above] {$250$} (3);
    \draw[->] (1) -- node[midway, below left] {$175$} (4);
    \draw[->] (3) -- (6);
    \draw[->] (4) -- (6);
    \draw[->] (5) -- node[midway, above right] {$432000$} (8);
    \draw[->] (6) -- node[midway, below] {$345600$} (8);
    
    \end{tikzpicture} 
    \caption{Example of a flow graph}
    \label{fig:example_flow}
\end{figure}
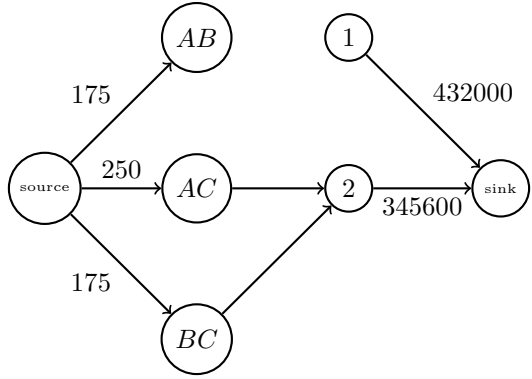

In this way, to calculate the $A_e^t$ parameters of our example, we can take the average session duration from Table \ref{table:example_data} per station over our 24-hour period. We assume for this example that all sessions are from the same day. This implies $110 \cdot 5 = 550$~kW for station 1 and $100 \cdot 3 + 10 \cdot 5 = 350$~kW for station 2. We sum all supplied kW within the same borough. Since station 1 and 2 are in different boroughs, the $\Lambda$ borough has 550~kW of supply and the $\Omega$ borough has 350~kW. Using the OD matrix of Table~\ref{table:od_matrix}, we can write two equations: $ 0.5 q_\Omega + 0.25 q_\Lambda = 350$ and $0.5 q_\Omega + 0.75 q_\Lambda = 550$. Solving this set of equations gives us $q_\Omega = 500$ and $q_\Lambda = 400$. This implies that we have 500~kW of demand in borough $\Omega$ and 400~kW in $\Lambda$. The resulting demand flow is as follow: between $\Omega$ and $\Omega$, $0.5 q_\Omega = 250$~kW, between $\Omega$ and $\Lambda$, $0.5 q_\Omega + 0.25 q_\Lambda = 350$~kW and between $\Lambda$ and $\Lambda$, $0.75 q_\Lambda = 300$~kW. We combine the flow from $\Omega$ to $\Lambda$ and from $\Lambda$ to $\Omega$, since we make the assumption that our flow is bidirectional. Finally, we split these flows between each point to obtain the flow for each OD pair: $AC$ has 250~kW, $AB$ has 175~kW and $BC$ has 175~kW. It is worth noting that in this example, we lost 300~kW of demand because we cannot represent the flow between $\Lambda$ and $\Lambda$ since we only have a single point in that borough. In practice, we generate a critical mass of points to guarantee at least 2 points in each borough. This gives us the final flow graph in Figure~\ref{fig:example_flow}.

We do not have a list of potential locations, so we use impossible demand to create candidate locations. Recall that impossible demand is the demand generated by an OD vertex with no edges to any station in the bipartite graph. This demand cannot be satisfied, and as such, considering a candidate location near it is likely a good possibility. To do so, we add two candidate locations for each impossible OD demand: one at the origin and one at the destination. Note that any OD vertex within the defined radius $R$ of a candidate location gets an outgoing arc to that location, including the OD vertex related with the impossible demand. In our example, the OD pair $AB$ does not have any valid station. As such, we would add two candidate locations: one at $A$ and one at $B$. This means both OD pairs $AC$ and $BC$ would also gain a new station. The kW/s given to the new station is based on an average across all stations of the same level. If we assume all stations in our example are of the same level then a new station would have $\approx$4.33~kW/s.

We do not have data about the cost of adding outlets or building new stations since the price can vary widely based on the location and electricity grid availability. As such, for our testing, we use arbitrary costs guided by reasonable considerations. We attribute to the addition of a level 2 outlet (to an existing station or a new one) a cost of 1. A level 3 outlet is twice that. Building a new level 2 station is 10 and a level 3 is 100. To account for these arbitrary costs, we run our experiments on multiple budgets ranging from 0 to 700 to perform a sensibility analysis of our MILP.

Having described our process for utilizing the data to generate the flow graphs and established the budget constraint, we now proceed to outline the various instances we create in accordance with the aforementioned procedure. We consider instances with $R \in \{400,500,600,700\}$, where the unit is meters, and with $W \in \{100,150,200,250,300\}$. For each combination of the $R$ and $W$ values, we generated 5 instances, where only the location of the $W$ random OD points differ. In each instance, the set $S_1$ contains the same 882 stations, while the set of candidate locations $S_2$ can vary depending on the OD pairs (recall the description above). For $W=100, 150, 200, 250$ and $300$, the instances have 4,950, 11,175, 19,900, 31,125 and 44,850 ODs (i.e., the cardinality of set $O$) and an average of 43.4, 65.8, 92.2, 117.4 and 141.8 candidate locations, respectively. In the next sections, we provide average results over the 5 generated instances for each $(R,W)$ pair. All the results are shown as a percentage of the total (average) kW charging demand, i.e., $\displaystyle\sum_{t \in T}\sum_{e \in L} A_e^t$.

\subsection{Assigning Demand to Stations}\label{subsec:AssignUsers}

Our first tests are meant to evaluate the impact of the radius $R$ and of the number of points $W$ on Program~\eqref{Program:linear}. The goal is to assess the sensitivity of satisfied and impossible demand to those parameters, analyze the service of the current infrastructure, as well as to identify a suitable value of $W$ that balances the model accuracy and computational efficiency when incorporating location and sizing decision. Note that larger values of $W$ allow a greater diversity of OD scenarios, but lead to larger optimization problems.

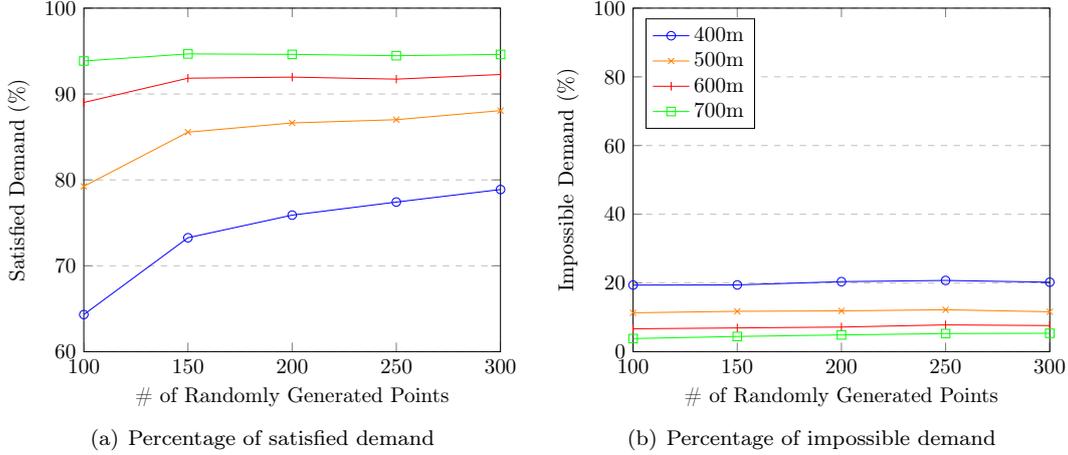
\begin{figure}  
\centering  

\subfigure[Percentage of satisfied demand]  
{  

\begin{tikzpicture}[scale=0.8]
\begin{axis}[
    xlabel={\# of Randomly Generated Points},
    ylabel={Satisfied Demand (\%)},
    xmin=100, xmax=300,
    ymin=60, ymax=100,
    xtick={100,150,200,250,300},
    ytick={0,20,40,60,70, 80, 90 ,100},
    ymajorgrids=true,
    grid style=dashed,
]

\addplot[
    color=blue,
    mark=o,
    ]
    coordinates {
(100,64.3210669008478
)(150,73.2636934263732
)(200,75.8974684076071
)(250,77.4191902536494
)(300,78.8788714017191
)
    };

\addplot[
    color=orange,
    mark=x,
    ]
    coordinates {
(100,79.248957391969
)(150,85.5556028193948
)(200,86.634296966989
)(250,87.0122887088231
)(300,88.0661644174032
)
    };

\addplot[
    color=red,
    mark=|,
    ]
    coordinates {
(100,89.013646547666
)(150,91.8404928621544
)(200,91.9673424261205
)(250,91.7306456166801
)(300,92.2724395419997
)
    };

\addplot[
    color=green,
    mark=square,
    ]
    coordinates {
(100,93.8430271933688
)(150,94.6670714729252
)(200,94.6138523474928
)(250,94.4581730093966
)(300,94.6070948727502
)
    };
    
\end{axis}
\end{tikzpicture}

}  
\subfigure[Percentage of impossible demand]  
{  
\label{subfig:bsingle}
\begin{tikzpicture}[scale=0.8] 
\begin{axis}[
    xlabel={\# of Randomly Generated Points},
    ylabel={Impossible Demand (\%)},
    xmin=100, xmax=300,
    ymin=0, ymax=100,
    xtick={100,150,200,250,300},
    ytick={0,20,40,60,80,100},
    legend pos=north west,
    ymajorgrids=true,
    grid style=dashed,
]

\addplot[
    color=blue,
    mark=o,
    ]
    coordinates {
(100,19.3850892847773
)(150,19.43174021775
)(200,20.3795027985346
)(250,20.7297338912572
)(300,20.2208683429178
)
    };
    \addlegendentry{400m}

\addplot[
    color=orange,
    mark=x,
    ]
    coordinates {
(100,11.3130812496239
)(150,11.74347756411
)(200,11.8753720499896
)(250,12.2022431082661
)(300,11.6108903477651
)
    };
    \addlegendentry{500m}

\addplot[
    color=red,
    mark=|,
    ]
    coordinates {
(100,6.65981511204519
)(150,6.90542242109996
)(200,7.16104523875205
)(250,7.79116552100678
)(300,7.5972045597276
)
    };
    \addlegendentry{600m}

\addplot[
    color=green,
    mark=square,
    ]
    coordinates {
(100,3.81231639687532
)(150,4.42160842070194
)(200,4.86675726657902
)(250,5.26205542339766
)(300,5.3451664134123
)
    };
    \addlegendentry{700m}
    
\end{axis}
\end{tikzpicture}
}
\caption{Single period results for Program~\eqref{Program:linear}}
\label{fig:24hourlinearresults}
\end{figure}

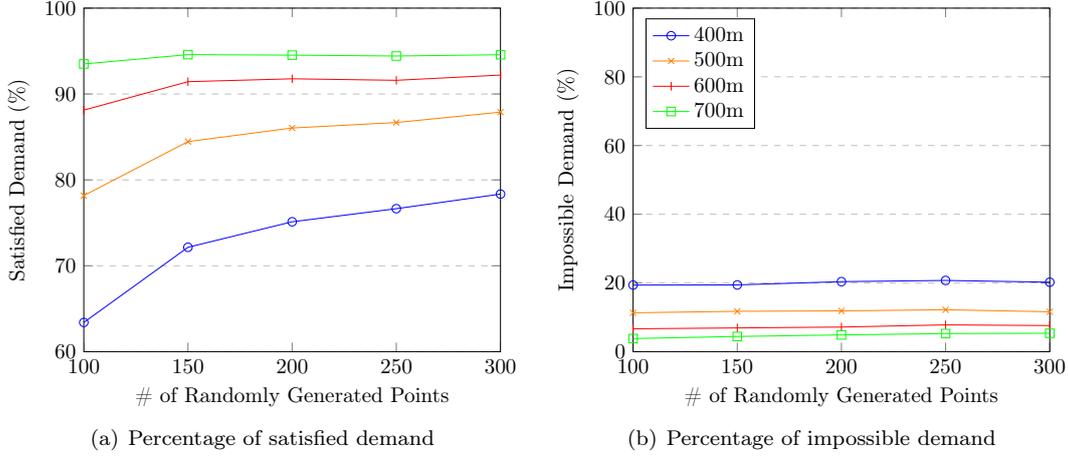
\begin{figure}  
\centering  

\subfigure[Percentage of satisfied demand]  
{  

\begin{tikzpicture}[scale=0.8] 
\begin{axis}[
    xlabel={\# of Randomly Generated Points},
    ylabel={Satisfied Demand (\%)},
    xmin=100, xmax=300,
    ymin=60, ymax=100,
    xtick={100,150,200,250,300},
    ytick={0,20,40,60,70, 80, 90 ,100},
    ymajorgrids=true,
    grid style=dashed,
]

\addplot[
    color=blue,
    mark=o,
    ]
    coordinates {
(100,63.4023688698981
)(150,72.1516413055341
)(200,75.1237106934892
)(250,76.6487643989779
)(300,78.3486925661528
)
    };

\addplot[
    color=orange,
    mark=x,
    ]
    coordinates {
(100,78.1747443130465
)(150,84.4618447193628
)(200,86.0438684064144
)(250,86.6814962353503
)(300,87.8901107366164
)
    };

\addplot[
    color=red,
    mark=|,
    ]
    coordinates {
(100,88.1261785485423
)(150,91.4364902931363
)(200,91.7704106691105
)(250,91.5954783670233
)(300,92.2045158630332
)
    };

\addplot[
    color=green,
    mark=square,
    ]
    coordinates {
(100,93.5000503487556
)(150,94.5829489948524
)(200,94.5242982624972
)(250,94.4235154138235
)(300,94.5714158413426
)
    };

\end{axis}
\end{tikzpicture}

}  
\subfigure[Percentage of impossible demand]  
{  

\begin{tikzpicture}[scale=0.8]
\begin{axis}[
    xlabel={\# of Randomly Generated Points},
    ylabel={Impossible Demand (\%)},
    xmin=100, xmax=300,
    ymin=0, ymax=100,
    xtick={100,150,200,250,300},
    ytick={0,20,40,60,80,100},
    legend pos=north west,
    ymajorgrids=true,
    grid style=dashed,
]

\addplot[
    color=blue,
    mark=o,
    ]
    coordinates {
(100,19.3850892847781
)(150,19.4317402177514
)(200,20.3795027985331
)(250,20.7297338912521
)(300,20.220868342902
)
    };
    \addlegendentry{400m}

\addplot[
    color=orange,
    mark=x,
    ]
    coordinates {
(100,11.3130812496245
)(150,11.743477564111
)(200,11.8753720499886
)(250,12.2022431082635
)(300,11.6108903477564
)
    };
    \addlegendentry{500m}

\addplot[
    color=red,
    mark=|,
    ]
    coordinates {
(100,6.65981511204552
)(150,6.90542242110056
)(200,7.16104523875151
)(250,7.79116552100513
)(300,7.59720455972198
)
    };
    \addlegendentry{600m}

\addplot[
    color=green,
    mark=square,
    ]
    coordinates {
(100,3.81231639687551
)(150,4.42160842070232
)(200,4.86675726657865
)(250,5.26205542339657
)(300,5.34516641340832
)
    };
    \addlegendentry{700m}
    
\end{axis}
\end{tikzpicture}
\label{subfig:b4blocks}
}

\caption{Multi-period results for Program~\eqref{Program:linear}}
\label{fig:6hourlinearresults}
\end{figure}

In Figure~\ref{fig:24hourlinearresults} (the y-axis of figure (a) begins at 60\% for better readability) and Figure~\ref{fig:6hourlinearresults}, we present our results for the satisfied and impossible demand\footnote{Recall that unsatisfied demand is distinct from both the satisfied and the impossible demand. Unsatisfied demand arises when the supply is insufficient.} separately for two cases: one where a single  period is considered, and the other where a day is discretized into 4 periods. In both cases, we observe that the satisfied demand increases as the number of points ($W$) increases, and then stabilizes. This is to be expected since the closer the number of points gets to the real number of EV users, the more accurately we depict the coverage of the territory by the existing infrastructure. Importantly, it appears that a relatively low value of $W$ suffices to capture the prevailing charging station coverage, which has a direct impact on the number of variables and constraints of Program~\eqref{Program:mip}, analyzed in the next section. Similarly, an increase in the radius $R$ results in an increase in satisfied demand. This can be attributed to the fact that a larger radius provides each OD with a greater number of station options to choose from. However, in contrast to the number of points $W$, opting for a larger radius distances us from reality, given that fewer individuals are expected to be willing to walk 700 meters compared to 400 meters. The larger radius can also be perceived as a compromise in service quality (coverage). Independently of that, any increase in the number of points or radius is limited by the amount of demand that can be satisfied. If the stations provide enough supply, the limiting factor becomes the impossible demand. It is worth noting that since our demand generation is closely related to the amount of supply (recall Section~\ref{subsec:MontrealCase}), it makes sense that we are able to satisfy most of the demand for these instances. On the flip side, increasing the radius reduces the impossible demand. This can be explained by the fact that a larger radius increases the covered area, meaning ODs are less likely to have no nearby stations. A less intuitive result is the fact that increasing the number of points does not affect the impossible demand. The reason behind this comes from the stochastic nature of our points' generation. If we have two instances with 100 points and 20\% of impossible demand each, if we assume that both instances do not have overlapping points which is statistically likely, then combining both instances into one leads to a 200 points instance with 20\% impossible demand. This only works if the 100 points instances are correctly predicting the total impossible demand of the network from the start.

The main difference between the two cases in these figures is in a slightly lower satisfied demand when we consider the four 6-hour periods. This is because the unsatisfied demand in each period cannot be satisfied in others, modeling an implicit constraint on when users are willing to charge. In the 24 hours instances, users can charge at any point of the day which is unrealistic. With four blocks of 6 hours, we can more easily reflect when users are looking for a station. However, it should be noted that adopting a more granular discretization would lead to more variables and constraints in our models. Moreover, this might also necessitate the incorporation of charging over consecutive periods (if time blocks are small). The impossible demand is exactly the same in both cases. This is because the impossible demand does not change within a day (we have the same OD pairs over the time horizon).

\begin{figure}[h]
    \centering
    \includegraphics[scale=0.5]{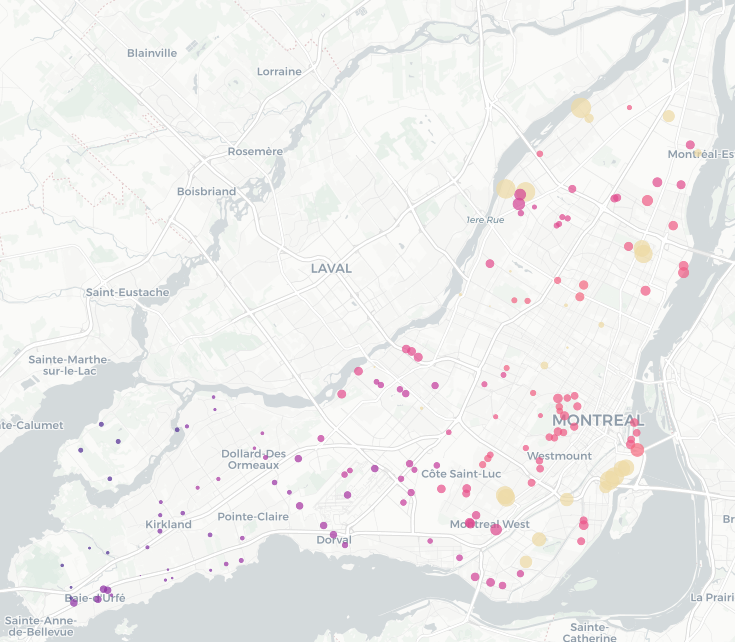}
    \caption{Map of the unsatisfied and impossible demand per point ($R=400, W=300$).}
    \medskip
    \footnotesize
     The points' size represents the unsatisfied demand and the colour represents the impossible demand (the larger and darker, the higher). Remark that points can have both unsatisfied and impossible demand as they are related to a set of OD pairs.
    \label{fig:unsatisfied_impossible_demand}
\end{figure}

Tables~\ref{table:percent_satisfied_level2} and \ref{table:percent_satisfied_level3} give a brief overview of the amount of demand satisfied at each period. The purpose of these results is to analyze their similarity with regards to the actual satisfied demand of Table~\ref{table:demandperperiod}. In other words, we aim to understand if the charging habits over the time periods are similar. In our results, we note that level 3 stations tend to mirror the satisfied demand of level 2 stations. This is because we do not model charging habits, i.e., users preferences. For instance,  during the period 0h-6h, level 3 satisfies more demand than in the other periods because in our instances (and session data) there is more demand over this period and there is no user preference making users to favour level 2. This could be improved to better match reality by either changing the capacity of level 3 stations during certain periods or introducing flow costs associated with the use of certain stations.  For the tables with the other radii, we refer to Appendix~\ref{app:AssignUsers}; the result trend is analogous. Figure~\ref{fig:unsatisfied_impossible_demand} shows a solved instance, where all points with unsatisfied and impossible demand are visible. The impossible demand becomes more prevalent as we go west, away from downtown. The closer we get to downtown, the impossible demand turns to unsatisfied demand. Near downtown, the unsatisfied demand is nonexistent with no points visible since all demand is satisfied.

\begin{table}[h!] \footnotesize
\caption {Percentage of demand per station for each time period (Figure~\ref{fig:demandperperiod})}
\label{table:demandperperiod}
\centering
\begin{tabular}{| c | c | c |} 
\hline
Period  & Level 2 & Level 3 \\
 \hline
  18h-24h & 23.64\%&23.91\%\\
 \hline
 12h-18h & 21.92\%& 28.35\%\\ 
 \hline
 6h-12h & 22.46\% & 24.82\%\\ 
 \hline
 0h-6h & 31.96\%& 22.91\% \\ 
 \hline
\end{tabular}
\end{table}

\begin{table}[h!] \footnotesize
\caption {Percentage of satisfied demand assigned to level 2 stations per time period with a 400 meter radius}
\label{table:percent_satisfied_level2}
\centering
\begin{tabular}{| c | c | c | c | c | c |} 
 \hline
 &\multicolumn{5}{c|}{\# of Randomly Generated Points} \\
Period  & 100 & 150 & 200 & 250 & 300 \\
 \hline
  18h-24h & 24.57\%&24.24\%& 23.90\%&23.66\%& 23.59\%\\
 \hline
 12h-18h & 23.44\%& 22.68\%& 22.20\%& 21.94\%& 21.74\%\\ 
 \hline
 6h-12h & 23.78\% & 23.08\%& 22.76\%&22.49\%& 22.32\%\\ 
 \hline
 0h-6h & 28.22\%& 30.00\%& 31.14\%& 31.91\% & 32.36\% \\ 
 \hline
\end{tabular}
\end{table}

\begin{table}[h!] \footnotesize
\caption {Percentage of satisfied demand assigned to level 3 stations per time period with a 400 meter radius}
\label{table:percent_satisfied_level3}
\centering
\begin{tabular}{| c | c | c | c | c | c |} 
 \hline
 &\multicolumn{5}{c|}{\# of Randomly Generated Points} \\
Period  & 100 & 150 & 200 & 250 & 300 \\
 \hline
  18h-24h &23.45\%&20.63\%& 22.92\%& 24.65\%&21.91\%\\
 \hline
 12h-18h & 20.27\%& 20.57\%& 19.98\%& 17.30\%& 18.97\%\\ 
 \hline
 6h-12h &21.77\% &21.46\%&20.55\%&20.96\%& 19.30\%\\ 
 \hline
 0h-6h & 34.51\% & 37.34\%& 36.55\%& 37.09\% & 39.81\%\\ 
 \hline
\end{tabular}
\end{table}

\newpage

\subsection{Adding and Expanding Stations}\label{subsec:AddingStations}

Our second set of tests is meant to evaluate the computational performance when solving our model, Program~\eqref{Program:mip}, for the installation of new stations and outlets. For these experiments, we exclusively focus on instances with $R=$ 400 meter radius, as it offers the greatest flexibility to users. In other words, this radius allows for a concentration of stations in close proximity to users. We take the instances with $W=$ 100, 150 and 200 random points to analyze the sensitivity of the optimal objective value (satisfied demand) to the number of points, while keeping the size of Program~\eqref{Program:mip} reasonable. As explained in Section~\ref{subsec:MontrealCase}, we use the impossible demand to identify a list of candidate locations $S_2$.

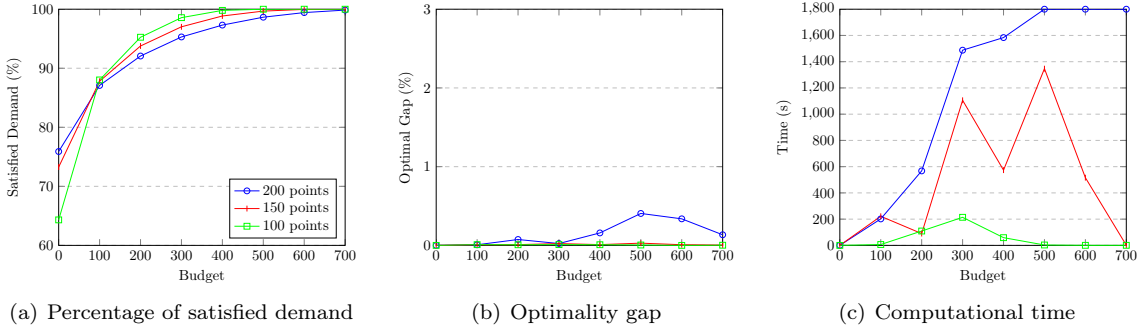
\begin{figure}[h]  
\centering  

\subfigure[Percentage of satisfied demand]  
{  

\begin{tikzpicture}[scale=0.55] 
\begin{axis}[
    xlabel={Budget},
    ylabel={Satisfied Demand (\%)},
    xmin=0, xmax=700,
    ymin=60, ymax=100,
    xtick={0,100,200,300,400,500,600,700},
    ytick={0,20,40,60,70, 80, 90 ,100},
    legend pos=south east,
    ymajorgrids=true,
    grid style=dashed,
]

\addplot[
    color=blue,
    mark=o,
    ]
    coordinates {
    (0,75.90)(100,87.09)(200,92.07)(300,95.30)(400,97.31)(500,98.6470175033446)(600,99.4273526994834)(700, 99.8486202661179)
    };
    \addlegendentry{200 points}

\addplot[
    color=red,
    mark=|,
    ]
    coordinates {
    (0,73.26)(100,87.80)(200,93.77)(300,97.02)(400,98.85)(500,99.6783929628579)(600,99.9643817905552)(700,99.9968731714031)
    };
    \addlegendentry{150 points}

\addplot[
    color=green,
    mark=square,
    ]
    coordinates {
    (0,64.32)(100,88.02)(200,95.25)(300,98.57)(400,99.80)(500,99.9950355939574)(600,100)(700,99.9979624917323)
    };
    \addlegendentry{100 points}
    
\end{axis}
\end{tikzpicture}

}  \hfill
\subfigure[Optimality gap]  
{  

\begin{tikzpicture}[scale=0.55] 
\begin{axis}[
    xlabel={Budget},
    ylabel={Optimal Gap (\%)},
    xmin=0, xmax=700,
    ymin=0, ymax=3.0,
    xtick={0,100,200,300,400,500,600,700},
    ytick={0.0,1.0,2.0,3.0},
    ymajorgrids=true,
    grid style=dashed,
]

\addplot[
    color=blue,
    mark=o,
    ]
    coordinates {
    (0,0)(100,0.00798644171233104)(200, 0.074706480500842)(300,0.0234938368331039)(400,0.158384259109591)(500,0.405543637607093)(600,0.337262268715364)(700,0.134634309254428)
    };

\addplot[
    color=red,
    mark=|,
    ]
    coordinates {
    (0,0)(100,0.00798429624669976)(200,0.00959465824107939)(300,0.0228405922971185)(400,0.00986150018505884)(500,0.0283332674347038)(600,0.00848076707927568)(700,0.00312709051519631)
    };

\addplot[
    color=green,
    mark=square,
    ]
    coordinates {
    (0,0)(100,0.00918325458991831)(200,0.00797577295051662)(300,0.0097549594101959)(400,0.00841460218782728)(500,0.00496472857403688)(600,0)(700,0.00203760042707856)
    };
    
\end{axis}
\end{tikzpicture}

}
\hfill
\subfigure[Computational time]  
{

\begin{tikzpicture}[scale=0.55]  
\begin{axis}[
    xlabel={Budget},
    ylabel={Time (s)},
    xmin=0, xmax=700,
    ymin=0, ymax=1800,
    xtick={0,100,200,300,400,500,600,700},
    ytick={0, 200, 400, 600, 800, 1000, 1200, 1400, 1600, 1800},
    ymajorgrids=true,
    grid style=dashed,
]

\addplot[
    color=blue,
    mark=o,
    ]
    coordinates {
    (0,0.347000000003027)(100,201.725199999986)(200, 567.659599999996)(300,1487.6724)(400,1583.6998)(500,1800)(600,1800)(700,1800)
    };

\addplot[
    color=red,
    mark=|,
    ]
    coordinates {
    (0,0.165399999974761)(100,221.828000000003)(200,90.0591999999946)(300,1105.84720000001)(400,573.406599999999)(500,1346.92799999999)(600,516.081199999992)(700,1.16259999999311)
    };

\addplot[
    color=green,
    mark=square,
    ]
    coordinates {
    (0,0.0689999999944121)(100,6.5564000000013)(200,109.462800000014)(300,213.975200000021)(400,58.8966000000015)(500,3.05)(600,0.402799999981653)(700,0.403200000000652)
    };

\end{axis}
\end{tikzpicture}

}

\caption{Single period results for Program~\eqref{Program:mip}}
\label{fig:single_period_stations}
\end{figure}

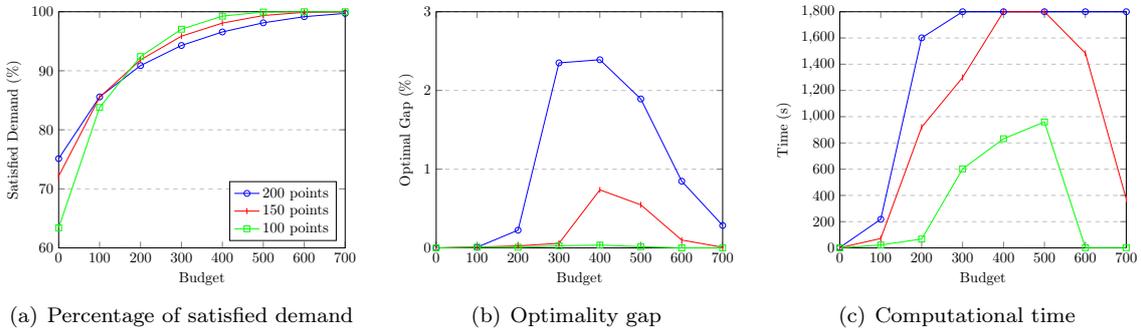
\begin{figure}  
\centering  

\subfigure[Percentage of satisfied demand]  
{  

\begin{tikzpicture}[scale=0.55] 
\begin{axis}[
    xlabel={Budget},
    ylabel={Satisfied Demand (\%)},
    xmin=0, xmax=700,
    ymin=60, ymax=100,
    xtick={0,100,200,300,400,500,600,700},
    ytick={0,20,40,60,70, 80, 90 ,100},
    legend pos=south east,
    ymajorgrids=true,
    grid style=dashed,
]

\addplot[
    color=blue,
    mark=o,
    ]
    coordinates {
    (0,75.1237106935014)(100,85.5712685451705)(200,90.8709015608073)(300,94.2837859321464)(400,96.5729613659013)(500,98.1215957450467)(600,99.1613246138241)(700,99.7157185677852)
    };
        \addlegendentry{200 points}

\addplot[
    color=red,
    mark=|,
    ]
    coordinates {
    (0,72.1516413055311)(100,85.5192968042404)(200,91.8435406876657)(300,95.8544983929412)(400,98.074768728318)(500,99.357354165842)(600,99.9000447214189)(700,99.9919640986224)
    };
        \addlegendentry{150 points}

\addplot[
    color=green,
    mark=square,
    ]
    coordinates {
    (0,63.4023688698984)(100,83.7879228678246)(200,92.4096199299392)(300,97.0362963224422)(400,99.2556471940438)(500,99.9450404563339)(600,100)(700,100)
    };
        \addlegendentry{100 points}
    
\end{axis}
\end{tikzpicture}

}  \hfill
\subfigure[Optimality gap]  
{  

\begin{tikzpicture}[scale=0.55]  
\begin{axis}[
    xlabel={Budget},
    ylabel={Optimal Gap (\%)},
    xmin=0, xmax=700,
    ymin=0, ymax=3,
    xtick={0,100,200,300,400,500,600,700},
    ytick={0.0,1.0,2.0,3.0},
    ymajorgrids=true,
    grid style=dashed,
]

\addplot[
    color=blue,
    mark=o,
    ]
    coordinates {
    (0,0)(100,0.00877109553045796)(200,0.225151362057662)(300,2.35032366707963)(400,2.39010950665751)(500,1.89206439165258)(600,0.847275834138188)(700,0.285283388064023)
    };

\addplot[
    color=red,
    mark=|,
    ]
    coordinates {
    (0,0)(100,0.00973758126936367)(200,0.0259906094520944)(300,0.0582307246550037)(400,0.73710339571469)(500,0.547401191859379)(600,0.100145369106223)(700,0.00803772857005012)
    };

\addplot[
    color=green,
    mark=square,
    ]
    coordinates {
    (0,0)(100,0.00900977756800282)(200,0.0082156515698106)(300,0.0268046088986389)(400,0.0360036262195449)(500,0.0154059180787144)(600,0)(700,0)
    };
    
\end{axis}
\end{tikzpicture}

}
\hfill
\subfigure[Computational time]  
{  

\begin{tikzpicture}[scale=0.55]  
\begin{axis}[
    xlabel={Budget},
    ylabel={Time (s)},
    xmin=0, xmax=700,
    ymin=0, ymax=1800,
    xtick={0,100,200,300,400,500,600,700},
    ytick={0, 200, 400, 600, 800, 1000, 1200, 1400, 1600, 1800},
    ymajorgrids=true,
    grid style=dashed,
]

\addplot[
    color=blue,
    mark=o,
    ]
    coordinates {
    (0,3.46219999999157)(100,218.23440000001)(200,1601.05300000001)(300,1800)(400,1800)(500,1800)(600,1800)(700,1800)
    };

\addplot[
    color=red,
    mark=|,
    ]
    coordinates {
    (0,1.37820000000065)(100,71.8344000000157)(200,920.906200000004)(300,1298.71599999999)(400,1800)(500,1800)(600,1483.7374)(700,370.765599999996)
    };

\addplot[
    color=green,
    mark=square,
    ]
    coordinates {
    (0,0.362800000002608)(100,22.228199999989)(200,67.728199999989)(300,602.053200000012)(400,831.875)(500,960.006199999992)(600,1.975)(700,2.09360000001034)
    };
    
\end{axis}
\end{tikzpicture}

}
\caption{Multi-period results for Program~\eqref{Program:mip}}
 \label{fig:multi_period_stations}
\end{figure}

Figures~\ref{fig:single_period_stations} and~\ref{fig:multi_period_stations} provide our results for the case with a single period and for the case with 4 periods. We can remark two trivial cases. The first is when the budget is 0, which simply results in the linear model since no station or outlet can be added. The second is when the budget is so large ($G \ge 700$) that we can add as many stations and outlets as necessary to satisfy all the demand, both unsatisfied and impossible. The interesting cases lie in the middle, where instances with different numbers of points $W$ have similar percentages. This can be explained by the fact that the impossible demand is relatively stable among instances (recall Figures~\ref{subfig:bsingle} and~\ref{subfig:b4blocks} for $R=$ 400~m) and, as such, adding a budget reduces impossible demand by similar amounts. We can observe that solving instances with an increasing number of points $W$, tens to lead to higher computational times and optimality gaps. Therefore, going beyond a $W$ value exceeding 200 points is anticipated to substantially increase computational times. The dips in the graphs are related to the high variance between instances (refer to the appendix \ref{app:AddingStations} for detailed results).

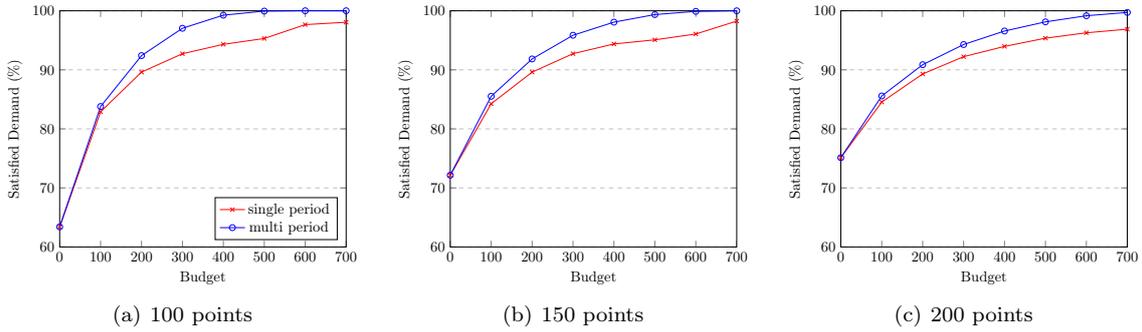
\begin{figure}
\subfigure[100 points]  
{  
\begin{tikzpicture} [scale=0.55] 
\begin{axis}[
    xlabel={Budget},
    ylabel={Satisfied Demand (\%)},
    xmin=0, xmax=700,
    ymin=60, ymax=100,
    xtick={0,100,200,300,400,500,600,700},
    ytick={0,20,40,60,70, 80, 90 ,100},
    legend pos=south east,
    ymajorgrids=true,
    grid style=dashed,
]
\addplot[
    color=red,
    mark=x,
    ]
    coordinates {
    (0,63.4023688698984)(100,82.89387591609666)(200,89.6225780490906)(300,92.7251154439565)(400,94.3166508121057)(500,95.3089930459817)(600,97.6582523756786)(700,98.0580263848792)
    };
    \addlegendentry{single period}
\addplot[
    color=blue,
    mark=o,
    ]
    coordinates {
    (0,63.4023688698984)(100,83.7879228678246)(200,92.4096199299392)(300,97.0362963224422)(400,99.2556471940438)(500,99.9450404563339)(600,100)(700,100)
    };
    \addlegendentry{multi period}

\end{axis}
\end{tikzpicture}
}\hfill
\subfigure[150 points]  
{  
\begin{tikzpicture} [scale=0.55] 
\begin{axis}[
    xlabel={Budget},
    ylabel={Satisfied Demand (\%)},
    xmin=0, xmax=700,
    ymin=60, ymax=100,
    xtick={0,100,200,300,400,500,600,700},
    ytick={0,20,40,60,70, 80, 90 ,100},
    ymajorgrids=true,
    grid style=dashed,
]
\addplot[
    color=red,
    mark=x,
    ]
    coordinates {
    (0,72.1516413055311)(100,84.2508254072011)(200,89.6099626485415)(300,92.7344907020773)(400,94.3790775794069)(500,95.0805825233)(600,96.0688029456677)(700,98.25311707608)
    };

\addplot[
    color=blue,
    mark=o,
    ]
    coordinates {
    (0,72.1516413055311)(100,85.5192968042404)(200,91.8435406876657)(300,95.8544983929412)(400,98.074768728318)(500,99.357354165842)(600,99.9000447214189)(700,99.9919640986224)
    };

\end{axis}
\end{tikzpicture}
}\hfill
\subfigure[200 points]  
{  
\begin{tikzpicture} [scale=0.55] 
\begin{axis}[
    xlabel={Budget},
    ylabel={Satisfied Demand (\%)},
    xmin=0, xmax=700,
    ymin=60, ymax=100,
    xtick={0,100,200,300,400,500,600,700},
    ytick={0,20,40,60,70, 80, 90 ,100},
    ymajorgrids=true,
    grid style=dashed,
]
\addplot[
    color=red,
    mark=x,
    ]
    coordinates {
    (0,75.1237106935014)(100,84.5414894840792)(200,89.3011865257274)(300,92.2276989142296)(400,93.9657252970515)(500,95.3598990532674)(600,96.2727010227566)(700,96.8848368360874)
    };

\addplot[
    color=blue,
    mark=o,
    ]
    coordinates {
    (0,75.1237106935014)(100,85.5712685451705)(200,90.8709015608073)(300,94.2837859321464)(400,96.5729613659013)(500,98.1215957450467)(600,99.1613246138241)(700,99.7157185677852)
    };

\end{axis}
\end{tikzpicture}
}
\caption{Single period results with multi-periods evaluation}
\label{fig:single_multi_period}
\end{figure}

In Figures~\ref{fig:single_period_stations} and~\ref{fig:multi_period_stations}, the main difference between the single and the multi-period cases is on the solving times and optimality gaps, with the multi-period one performing worst on both metrics. This is not surprising as the multi-period instances result in larger mixed-integer programs. With respect to the percentage of satisfied demand, the single and multi-period instances seems similar. In fact, the satisfied demand in the multi-period case is overall lower by at most 4.23\% and on average 1\% less. Although these differences are small, it is important to note that this comparison is not entirely fair, given that the instances are fundamentally distinct due to their utilization of different time discretizations. For this reason,  we evaluate the location and sizing decisions of solutions derived from single-period instances within the more realistic multi-period program, resulting in Figure \ref{fig:single_multi_period}.  We observe a discrepancy of up to 5.28\%  and an average of 2.75\% less satisfied demand compared to the multi-period solution (Figure~\ref{fig:multi_period_stations}). This shows that accounting for the time component of the charging location and sizing problem is meaningful. 

\begin{figure}[h]
    \centering
    \includegraphics[scale=0.5]{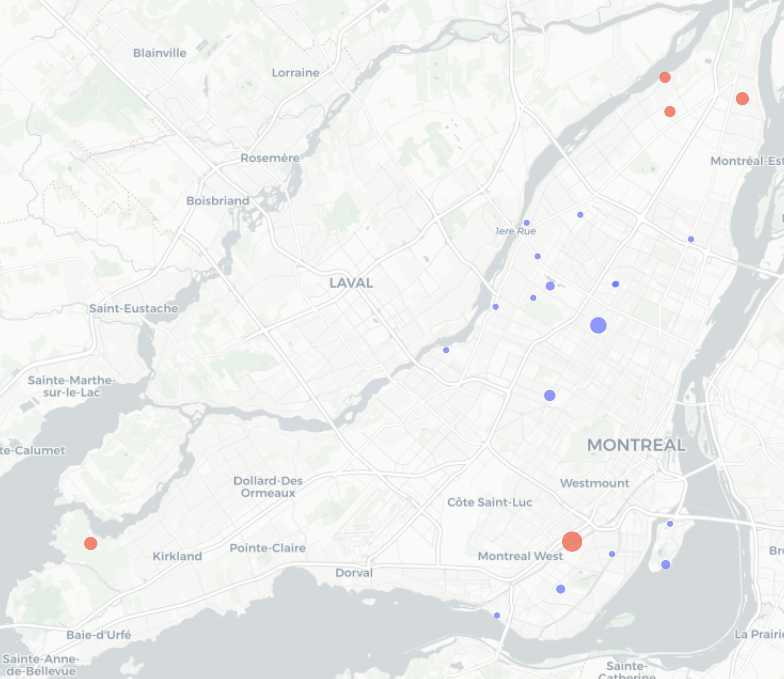}
    \caption{Map of Montreal with only new stations and new outlets (R=400, W=200, G=100)}
    \label{fig:added_stations}
    \medskip
    \footnotesize
    The blue points represents existing level 2 stations and the orange ones are new level 2 stations. The size of the points corresponds to the number of outlets added. No level 3 stations are expanded or opened.
\end{figure}

While analyzing the solutions to our instances, we observed that for the Montreal case, our model tends to prioritize the addition of level 2 outlets to existing stations and the addition of level 2 stations, as we increase the budget. In fact, among the solved instances, no solution included the addition of level 3 outlets or stations. This can be explained by the fact that the existing infrastructure is mostly sufficient to sustain the already existing demand. As such, the majority of the budget is spent providing a sparse amount of supply in undersupplied areas. Based on our cost parameters, level 3 stations are simply too expensive for this specific use case. Figure \ref{fig:added_stations} is an example of this behavior.

\section{Conclusion}\label{sec:conclusion}

In this work, we propose a maximum flow-based formulation for modelling the problem of locating and sizing EV charging stations in an intracity context. By solving our mathematical programming model, we are able to optimize the placement of new stations as well as to evaluate their expected usage (i.e., the charging demand they satisfy), within a densely populated city. Performing such an evaluation can be computationally expensive due to the large number of users, existing stations and candidate locations. However, our key contribution, the transformation of the charging demand assignment to stations into a maximum flow problem, allows for the scalability of our approach. To the best of our knowledge, we present the first approach capable of handling realistically-sized instances, multi-period instances with an exact algorithm. This improves upon single-period, exact methodologies (e.g., \cite{Shahraki2015}, \cite{Yang2017}) and multi-period, heuristic methodologies (e.g., \cite{Zhang2017}, \cite{Anjos2020}) from the existing literature.

Our mathematical programming models are run based on real-world data. Concretely, the linear model can evaluate the quality of service provided by an already existing network of stations in terms of satisfied demand. Moreover, it can be used to identify impossible demand and hence, locations where the installation of new stations would be desirable. The mixed-integer model can take a list of candidate locations and return the most promising ones or increase the capacity of the already existing infrastructure given a budget. We also demonstrate the impact of discretizing a 24-hour period to account for highs and lows in demand over a day.

For the application of our approach to the case study of the island of Montreal, the data needed was the charging session records of all the existing infrastructure and the daily OD travels across all boroughs. We show that the current infrastructure would not require a large increase in supply to satisfy all the demand at the moment. However, despite providing enough supply, the network does not provide uniform quality of service across the island, resulting in certain regions having limited access to public charging facilities.

In practice, our methodology should be especially useful for infrastructure owners to identify limitations in their provision of charging, namely, regions with impossible and unsatisfied demand. The direct use of our optimal expansion decisions must be cautiously analyzed for each specific application as simplifications were made for sake of tractability.

Further research could focus on the integration of our maximum flow model into  EV charging stations placement problems using other objectives such as cost minimization. This would result in a bilevel program with a maximum flow problem at the lower level. Another important aspect would be the integration of power grid constraints, which could restrict candidate locations. Concerning the modeling of the assignment of the demand to stations, an aspect for future consideration is to account for user preferences over stations (or locations), instead of assuming that the demand is effectively spread (e.g., through a real-time app). This could be potentially done by assigning weights to the arcs of our flow network. Another line on research could be to further explore the best way to discretize a day to properly reflect reality.  On the same topic, although we consider a time horizon for our flow model, we do not allow for flow to travel between consecutive periods, which could occur in practice.  Finally, expanding our approach to handle both the intracity and intercity cases would allow for a more complete modelling of the optimal location and sizing of EV charging stations problem.

 %
 %



\section*{Acknowledgments}
We are grateful to Jean-Luc Dupr\'e from the \emph{Direction Mobilit\'e} of \emph{Hydro-Qu\'ebec} for sharing his knowledge about our case study, and to Steven Lamontagne for his insights about the project.
    
This research was supported by Hydro-Québec, NSERC Collaborative Research and Development Grant CRDPJ 536757 - 19, and the FRQ-IVADO Research Chair in Data Science for Combinatorial Game Theory.

All the map figures presented in this paper were created using \texttt{geojson.io}\footnote{https://github.com/mapbox/geojson.io} and  the open source graphing library for Python \texttt{plotly}\footnote{https://plotly.com/python/}.





\bibliographystyle{unsrt}
\bibliography{ref}



\newpage

\appendix


\section{Futher Computational Results}

\subsection{Detailed Results of Section~\ref{subsec:AssignUsers}} \label{app:AssignUsers}

\begin{figure}[h]
    \centering
    \subfigure[0h-6h]  
{  
    \includegraphics[scale=0.3]{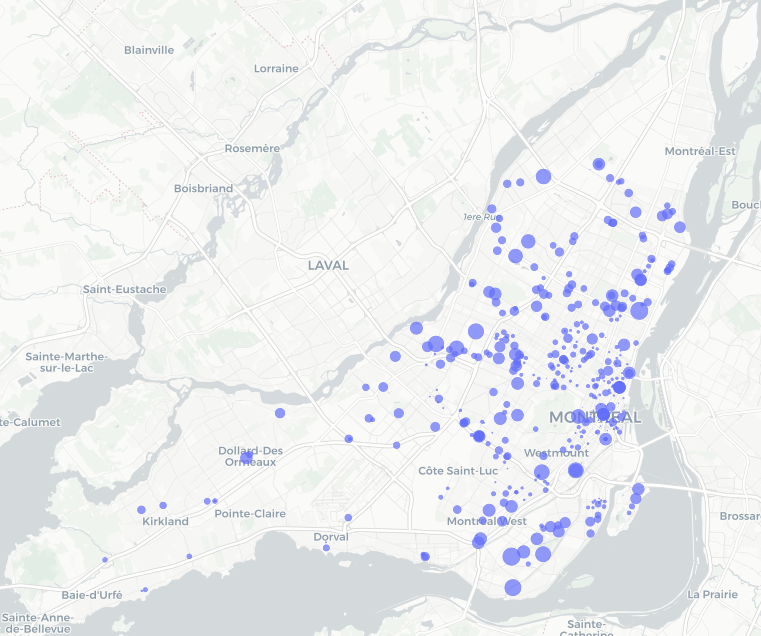}
}
\subfigure[6h-12h]  
{  
    \includegraphics[scale=0.3]{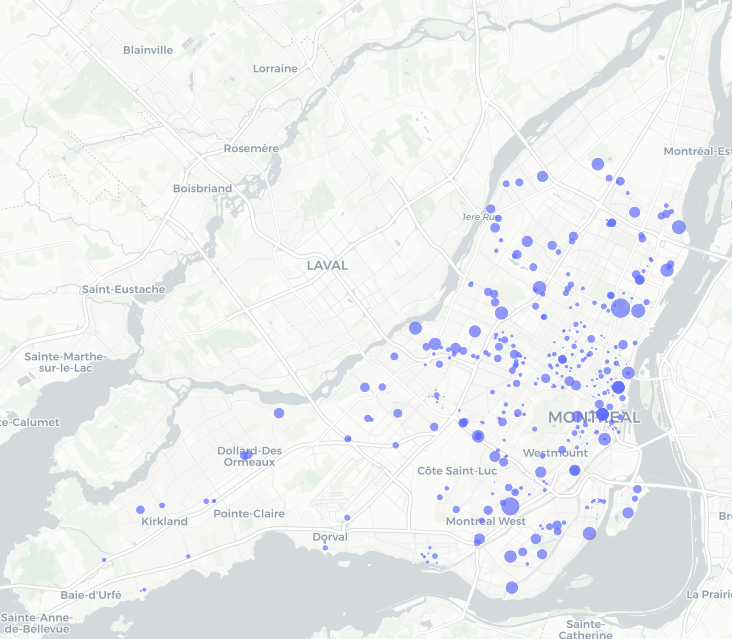}
}
\subfigure[12h-18h]  
{  
    \includegraphics[scale=0.3]{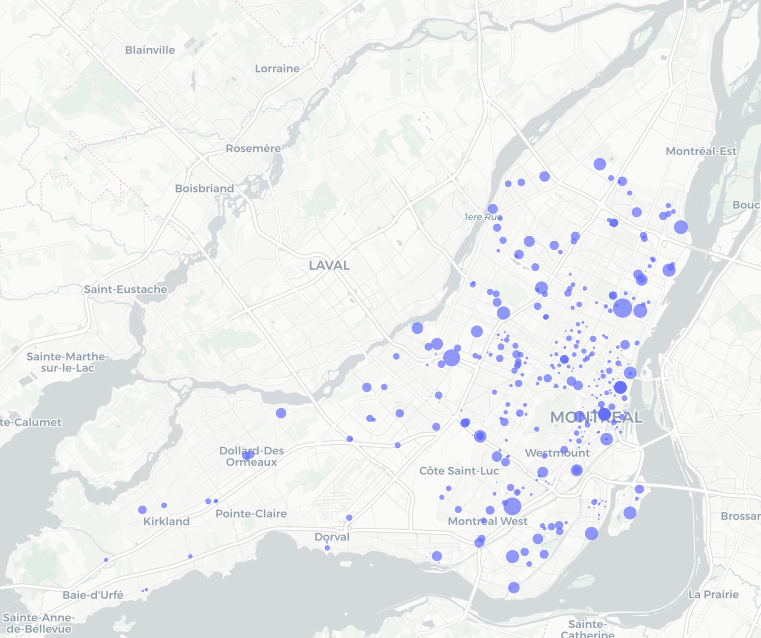}
}
\subfigure[18h-24h]  
{  
    \includegraphics[scale=0.3]{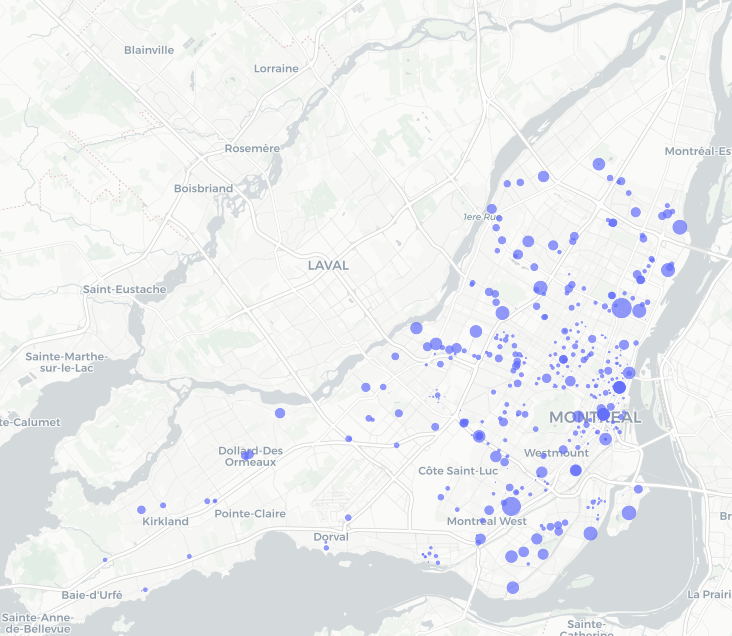}
}
    \caption{Example of satisfied demand over time periods (R=700, W=300)}
    \label{fig:heatmap}
\end{figure}

Figure~\ref{fig:heatmap} provides an illustration on the variation of the satisfied demand over each block of 6 hours. Tables~\ref{table:A1} to~\ref{table:A6} provide further details on the average satisfied demand for each value of $R$.

\begin{table}[h!] \footnotesize
\caption {Percentage of satisfied demand assigned to level 2 stations per time period with a 500 meter radius}
\begin{tabular}{| c | c | c | c | c | c |} 
 \hline
 &\multicolumn{5}{c|}{\# of Randomly Generated Points} \\
Period  & 100 & 150 & 200 & 250 & 300 \\
 \hline
  18h-24h &24.17\%&23.90\%& 23.59\%& 23.44\%&23.41\%\\
 \hline
 12h-18h & 22.60\%&22.04\%& 21.76\%& 21.64\%& 21.54\%\\ 
 \hline
 6h-12h &23.13\% &22.54\%&22.32\%&22.19\%& 22.16\%\\ 
 \hline
 0h-6h & 30.10\% & 31.52\%& 32.34\%& 32.72\% & 32.89\%\\ 
 \hline
\end{tabular}\label{table:A1}
\end{table}

\begin{table}[h!] \footnotesize
\caption {Percentage of satisfied demand assigned to level 3 stations per time period with a 500 meter radius}
\begin{tabular}{| c | c | c | c | c | c |} 
 \hline
 &\multicolumn{5}{c|}{\# of Randomly Generated Points} \\
Period  & 100 & 150 & 200 & 250 & 300 \\
 \hline
  18h-24h &23.23\%&21.46\%& 22.17\%& 23.71\%&22.20\%\\
 \hline
 12h-18h & 21.24\%&20.31\%& 20.05\%& 18.73\%& 19.53\%\\ 
 \hline
 6h-12h &21.34\% &21.92\%&20.79\%&19.89\%& 20.36\%\\ 
 \hline
 0h-6h & 34.19\% & 36.31\%& 36.98\%& 37.67\% & 37.91\%\\ 
 \hline
\end{tabular}
\end{table}

\begin{table}[h!] \footnotesize
\caption {Percentage of satisfied demand assigned to level 2 stations per time period with a 600 meter radius}
\begin{tabular}{| c | c | c | c | c | c |} 
 \hline
 &\multicolumn{5}{c|}{\# of Randomly Generated Points} \\
Period  & 100 & 150 & 200 & 250 & 300 \\
 \hline
  18h-24h &23.68\%&23.39\%& 23.44\%&23.36\%&23.25\%\\
 \hline
 12h-18h & 22.02\%&21.52\%& 21.54\%& 21.62\%& 21.43\%\\ 
 \hline
 6h-12h &22.60\% &22.10\%&22.21\%&22.17\%& 22.02\%\\ 
 \hline
 0h-6h & 31.70\%& 32.99\%& 32.82\%& 32.85\% & 33.30\%\\ 
 \hline
\end{tabular}
\end{table}

\begin{table}[h!] \footnotesize
\caption {Percentage of satisfied demand assigned to level 3 stations per time period with a 600 meter radius}
\begin{tabular}{| c | c | c | c | c | c |} 
 \hline
 &\multicolumn{5}{c|}{\# of Randomly Generated Points} \\
Period  & 100 & 150 & 200 & 250 & 300 \\
 \hline
  18h-24h &24.14\%&23.61\%& 22.71\%& 23.02\%&23.81\%\\
 \hline
 12h-18h & 21.55\%&22.09\%& 21.20\%& 19.58\%& 21.10\%\\ 
 \hline
 6h-12h &21.67\% &22.71\%&21.15\%&20.68\%& 21.80\%\\ 
 \hline
 0h-6h & 32.64\% & 31.60\%& 34.94\%& 36.71\% & 33.29\%\\ 
 \hline
\end{tabular}
\end{table}

\begin{table}[h!] \footnotesize
\caption {Percentage of satisfied demand assigned to level 2 stations per time period with a 700 meter radius}
\begin{tabular}{| c | c | c | c | c | c |} 
 \hline
 &\multicolumn{5}{c|}{\# of Randomly Generated Points} \\
Period  & 100 & 150 & 200 & 250 & 300 \\
 \hline
  18h-24h &23.27\%&23.35\%& 23.27\%&23.34\%&23.21\%\\
 \hline
 12h-18h & 21.45\%&21.45\%& 21.44\%& 21.58\%& 21.52\%\\ 
 \hline
 6h-12h &22.06\%&22.05\%&22.10\%&22.11\%& 22.13\%\\ 
 \hline
 0h-6h & 33.22\%& 33.16\%& 33.19\%& 32.97\% &33.13\%\\ 
 \hline
\end{tabular}
\end{table}

\begin{table}[h!] \footnotesize
\caption {Percentage of satisfied demand assigned to level 3 stations per time period with a 700 meter radius}
\begin{tabular}{| c | c | c | c | c | c |} 
 \hline
 &\multicolumn{5}{c|}{\# of Randomly Generated Points} \\
Period  & 100 & 150 & 200 & 250 & 300 \\
 \hline
  18h-24h &23.99\%&22.76\%& 23.70\%& 22.87\%&23.95\%\\
 \hline
 12h-18h & 21.88\%&21.19\%& 21.55\%& 19.99\%& 20.41\%\\ 
 \hline
 6h-12h &22.40\%&21.77\%&21.77\%&21.44\%& 21.12\%\\ 
 \hline
 0h-6h & 31.73\%& 34.28\%& 32.98\%& 35.69\% & 34.52\%\\ 
 \hline
\end{tabular}\label{table:A6}
\end{table}

\newpage

\subsection{Detailed Results for Section~\ref{subsec:AddingStations}} \label{app:AddingStations}

Tables~\ref{table:A7} to~\ref{table:A12} provide detailed results with respect to the satisfied demand, solving times and optimality gap for each of our instances. 

\begin{table}[h!] \footnotesize
\caption {Satisfied demand for Program~\eqref{Program:mip} over a single period}
\begin{tabular}{| c | c | c | c | c | c | c | c | c |} 
 \hline
  & \multicolumn{8}{|c|}{Budget} \\
  & 0 & 100 & 200 & 300 & 400 & 500 & 600 & 700\\
 \hline
 \multicolumn{9}{|c|}{100 Points} \\
 \hline
 Instance 1 & 62.38\% & 89.14\% & 95.48\%& 98.85\% & 99.86\%& 99.99\%& 100.00\%&100.00\%\\
 \hline
 Instance 2 & 67.84\% &88.43\% & 95.79\% & 98.73\% &99.90\% & 99.99\%& 100.00\% & 100.00\% \\
 \hline
 Instance 3 & 57.72\% & 84.21\% & 93.54\% & 97.88\% & 99.64\% & 99.99\%& 100.00\% & 100.00\% \\
 \hline
 Instance 4 & 67.81\%& 92.85\% & 97.90\% &99.79\% & 100.00\% & 100.00\% & 100.00\% & 99.99\% \\
 \hline
 Instance 5 & 65.85\%& 85.46\% & 93.54\% & 97.62\% & 99.60\% & 100.00\% & 100.00\% & 100.00\% \\
 \hline
 \multicolumn{9}{|c|}{150 Points} \\
 \hline
 Instance 1 & 76.49\% & 90.04\% &95.07\% & 97.73\%& 99.24\%& 99.87\% & 99.99\% &100.00\%\\
 \hline
 Instance 2 & 69.97\% & 86.27\% &93.51\% & 97.23\% & 98.90\%& 99.75\% & 99.99\%& 100.00\%\\
 \hline
 Instance 3 & 67.54\% &83.65\%& 91.04\%& 95.15\% & 97.85\% & 99.24\% & 99.88\% & 99.99\% \\
 \hline
 Instance 4 & 76.47\% & 91.74\% & 95.81\% & 98.11\% & 99.44\% & 99.91\% & 100.00\% &100.00\% \\
 \hline
 Instance 5 & 75.85\% & 87.30\% & 93.43\% & 96.86\% & 98.81\% & 99.63\% & 99.96\% & 99.99\% \\
 \hline
 \multicolumn{9}{|c|}{200 Points} \\
 \hline
 Instance 1 & 77.00\% & 88.84\% & 92.87\% & 95.91\% & 97.61\% & 98.84\% & 99.53\% & 99.92\% \\
 \hline
 Instance 2 & 74.86\% & 86.25\% & 91.79\% & 95.35\% & 97.42\% & 98.66\% & 99.47\%& 99.88\% \\
 \hline
 Instance 3 & 71.38\% & 83.16\%& 89.46\% & 93.48\% & 96.01\%& 97.92\% & 99.02\% & 99.67\% \\
 \hline
 Instance 4 & 79.68\% & 91.13\% & 94.77\% &96.90\%& 98.45\% & 99.27\% & 99.79\%& 99.97\% \\
 \hline
 Instance 5 & 76.56\% & 86.06\% & 91.46\% & 94.88\% & 97.07\% & 98.54\% & 99.32\% & 99.81\% \\
 \hline
\end{tabular}\label{table:A7}
\end{table}

\begin{table}[h!] \footnotesize
\caption {Computing time (s) for Program~\eqref{Program:mip} over a single period}
\begin{tabular}{| c | c | c | c | c | c | c | c | c |} 
 \hline
  & \multicolumn{8}{|c|}{Budget} \\
  & 0 & 100 & 200 & 300 & 400 & 500 & 600 & 700\\
 \hline
 \multicolumn{9}{|c|}{100 Points} \\
 \hline
 Instance 1 & 0.078 & 7.156 & 447.078 & 80.485 & 95.047 & 0.546 & 0.437 & 0.469 \\
 \hline
 Instance 2 & 0.063 & 7.625 & 10.079 & 27.375 & 19.578 & 1.016 & 0.5 & 0.406 \\
 \hline
 Instance 3 & 0.078 & 6.672 & 16.609 & 29.156 & 91.062 & 7.297 & 0.343 & 0.375 \\
 \hline
 Instance 4 & 0.063 & 5.282 & 55.266 & 14.563 & 0.281 & 0.328 & 0.328 & 0.422 \\
 \hline
 Instance 5 & 0.063 & 6.047 & 18.282 & 926.094 & 88.515 & 6.063 & 0.406 & 0.344 \\
 \hline
 \multicolumn{9}{|c|}{150 Points} \\
 \hline
 Instance 1 & 0.14 & 413.594 & 103.828 & 1800 & 861.797 & 1536.781 & 23.718 & 0.735 \\
 \hline
 Instance 2 & 0.172 & 33.578 & 133.328 & 64.75 & 657.157 & 655.625 & 830.922 & 0.906 \\
 \hline
 Instance 3 & 0.172 & 418.547 & 42.125 & 1800 & 587.137 & 1800 & 1031.094 & 2.453 \\
 \hline
 Instance 4 & 0.171 & 14.468 & 78.828 & 959.875 & 503.063 & 1800 & 1.578 & 0.75 \\
 \hline
 Instance 5 & 0.172 &  228.953 & 92.187 & 312.156 & 258.016 & 942.156 & 693.094 & 0.969 \\
 \hline
 \multicolumn{9}{|c|}{200 Points} \\
 \hline
 Instance 1 & 0.328 & 254.719 & 1800 & 1800 & 1800 & 1800 & 1800 & 1800 \\
 \hline
 Instance 2 & 0.282 & 11.719 & 413.36 & 1224.313 & 1428.156 & 1800 & 1800 & 1800 \\
 \hline
 Instance 3 & 0.421 & 665.766 & 468.86 & 1800 & 1800 & 1800 & 1800 & 1800 \\
 \hline
 Instance 4 & 0.36 & 12.235 & 61.734 & 1800 & 1090 & 1800 & 1800 & 1800 \\
 \hline
 Instance 5 & 0.344 & 64.187 & 94.297 & 813.891 & 1800 & 1800 & 1800 & 1800 \\
 \hline
\end{tabular}
\end{table}

\begin{table}[h!] \footnotesize
\caption {Optimal gap for Program~\eqref{Program:mip} over a single period}
\begin{tabular}{| c | c | c | c | c | c | c | c | c |} 
 \hline
  & \multicolumn{8}{|c|}{Budget} \\
  & 0 & 100 & 200 & 300 & 400 & 500 & 600 & 700\\
 \hline
 \multicolumn{9}{|c|}{100 Points} \\
 \hline
 Instance 1 & 0.00\% & 0.01\% & 0.01\% & 0.01\% & 0.01\% & 0.00\% &0.00\% & 0.00\% \\
 \hline
 Instance 2 &0.00\% & 0.01\% & 0.01\% & 0.01\% & 0.01\% & 0.01\% & 0.00\%  & 0.00\%  \\
 \hline
 Instance 3 & 0.00\%  & 0.01\% & 0.01\% &0.01\%& 0.01\% & 0.01\% & 0.00\% & 0.00\%\\
 \hline
 Instance 4 & 0.00\% & 0.01\% & 0.01\%& 0.01\% & 0.00\% & 0.00\% &0.00\% & 0.01\% \\
 \hline
 Instance 5 & 0.00\%  & 0.01\% & 0.00\% & 0.01\%  & 0.01\%  & 0.00\%& 0.00\% & 0.00\% \\
 \hline
 \multicolumn{9}{|c|}{150 Points} \\
 \hline
 Instance 1 & 0.00\% & 0.01\% & 0.01\% & 0.02\% & 0.01\% & 0.01\% & 0.01\% & 0.00\% \\
 \hline
 Instance 2 & 0.00\% & 0.00\% & 0.01\% & 0.01\% & 0.01\% &0.01\% & 0.01\% & 0.00\% \\
 \hline
 Instance 3 & 0.00\% & 0.01\% & 0.01\% & 0.07\% & 0.01\% & 0.10\% & 0.01\% & 0.01\% \\
 \hline
 Instance 4 & 0.00\% &  0.01\% & 0.01\%  &  0.01\%  &  0.01\%  &  0.01\% & 0.00\% & 0.00\% \\
 \hline
 Instance 5 & 0.00\% & 0.01\%& 0.01\% & 0.01\% & 0.01\% & 0.01\% & 0.01\% &0.01\% \\
 \hline
 \multicolumn{9}{|c|}{200 Points} \\
 \hline
 Instance 1 & 0.00\% & 0.01\% & 0.33\% & 0.02\% & 0.45\% &0.84\% & 0.36\% &0.08\% \\
 \hline
 Instance 2 & 0.00\% & 0.00\% & 0.01\% &0.01\% & 0.01\% & 0.48\%& 0.24\% & 0.09\% \\
 \hline
 Instance 3 & 0.00\% & 0.01\% & 0.01\% & 0.01\% & 0.23\%& 0.19\% & 0.58\% & 0.33\% \\
 \hline
 Instance 4 & 0.00\% & 0.01\% & 0.01\% & 0.06\% & 0.01\% & 0.42\% & 0.13\% & 0.02\% \\
 \hline
 Instance 5 & 0.00\% & 0.01\% & 0.01\% & 0.01\%& 0.09\%& 0.09\% & 0.38\% & 0.14\% \\
 \hline
\end{tabular}
\end{table}

\begin{table}[h!] \footnotesize
\caption {Satisfied demand for Program~\eqref{Program:mip} over multiple periods}
\begin{tabular}{| c | c | c | c | c | c | c | c | c |} 
 \hline
  & \multicolumn{8}{|c|}{Budget} \\
  & 0 & 100 & 200 & 300 & 400 & 500 & 600 & 700\\
 \hline
 \multicolumn{9}{|c|}{100 Points} \\
 \hline
 Instance 1 & 62.28\% & 84.66\% & 92.94\% & 97.24\% & 99.36\%& 99.97\% & 100.00\% & 100.00\% \\
 \hline
 Instance 2 & 66.20\% & 84.82\% & 93.56\% &97.35\% & 99.56\% & 99.98\% & 100.00\% & 100.00\% \\
 \hline
 Instance 3 & 57.30\% & 79.99\% & 89.82\% & 95.76\% & 98.73\% & 99.91\% & 100.00\%& 100.00\% \\
 \hline
 Instance 4 & 67.32\% & 89.24\% & 96.26\% &99.12\% & 99.97\% & 100.00\%& 100.00\% & 100.00\% \\
 \hline
 Instance 5 &63.91\% & 80.23\%& 89.48\%& 95.71\%& 98.66\%& 99.87\% & 100.00\% & 100.00\% \\
 \hline
 \multicolumn{9}{|c|}{150 Points} \\
 \hline
 Instance 1 & 74.87\% & 87.65\% & 93.26\% & 96.79\% & 98.59\% & 99.60\%& 99.97\% & 100.00\%\\
 \hline
 Instance 2 & 69.85\% & 84.17\% & 91.67\% & 96.09\% & 98.19\% & 99.49\% &99.93\%&100.00\% \\
 \hline Instance 3 & 66.03\% & 80.83\% & 88.18\% & 93.23\% & 96.55\%& 98.57\% & 99.63\% &99.97\%\\
 \hline
 Instance 4 & 75.62\% & 89.24\% & 94.38\% & 97.23\% & 98.87\% & 99.73\% & 99.99\% & 100.00\% \\
 \hline
 Instance 5 &74.39\% & 85.71\%& 91.73\% & 95.92\% & 98.17\% & 99.40\% & 99.87\% & 99.99\% \\
 \hline
 \multicolumn{9}{|c|}{200 Points} \\
 \hline
 Instance 1 & 76.12\% & 87.07\% & 91.73\% & 94.82\%& 97.00\% & 98.13\% & 99.34\% & 99.78\% \\
 \hline
 Instance 2 & 74.03\% &84.38\%& 90.70\% & 94.28\% & 96.50\%& 98.13\% & 99.17\% &99.76\% \\
 \hline
 Instance 3 & 70.63\% &81.66\% & 87.76\% & 91.96\% & 95.03\% & 97.21\% & 98.50\%&99.46\% \\
 \hline
 Instance 4 & 78.51\% & 89.43\%& 93.77\% & 96.27\% & 97.88\% &99.00\% & 99.67\% & 99.87\% \\
 \hline
 Instance 5 & 76.33\% &85.31\% &90.39\%& 94.09\%& 96.45\%& 98.14\% & 99.12\% & 99.72\% \\
 \hline
\end{tabular}
\end{table}

\begin{table}[h!] \footnotesize
\caption {Computing time (s) for Program~\eqref{Program:mip} over multiple periods}
\begin{tabular}{| c | c | c | c | c | c | c | c | c |} 
 \hline
  & \multicolumn{8}{|c|}{Budget} \\
  & 0 & 100 & 200 & 300 & 400 & 500 & 600 & 700\\
 \hline
 \multicolumn{9}{|c|}{100 Points} \\
 \hline
 Instance 1 & 0.313 & 9.657 & 34.235 & 720.25 & 1037.297 & 1308.797 & 1.922 & 1.781 \\
 \hline
 Instance 2 & 0.344 & 20.422 & 52.828 & 1800 & 222.047 & 946.719 & 1.844 & 1.906 \\
 \hline
 Instance 3 & 0.422 & 10.203 & 41.359 & 210.234 & 833.25 & 743.031 & 2.234 & 2.25 \\
 \hline
 Instance 4 & 0.375 & 19.437 & 24.532 & 97.797 & 266.703 & 1.406 & 1.625 & 2.188 \\
 \hline
 Instance 5 & 0.36 & 51.422 & 185.687 & 181.938 & 1800 & 1800 & 2.25 & 2.343 \\
 \hline
 \multicolumn{9}{|c|}{150 Points} \\
 \hline
 Instance 1 & 0.922 & 152.844 & 1006.859 & 1421.375 & 1800 & 1800 & 1800 & 9.109 \\
 \hline
 Instance 2 & 1.39 & 55.203 & 455.813 & 802.36 & 1800 & 1800 & 1800 & 29.891 \\
 \hline
 Instance 3 & 1.047 & 43.203 & 1116.75 & 1800 & 1800 & 1800 & 1800 & 1800 \\
 \hline
 Instance 4 & 1.969 & 34.469 & 224.922 & 669.594 & 1800 & 1800 & 218.344 & 4.781 \\
 \hline
 Instance 5 & 1.563 & 73.453 & 1800 & 1800 & 1800 & 1800 & 1800 & 9.969 \\
 \hline
 \multicolumn{9}{|c|}{200 Points} \\
 \hline
 Instance 1 & 2.156 & 496.344 & 1800 & 1800 & 1800 & 1800 & 1800 & 1800 \\
 \hline
 Instance 2 & 2.046 & 158.906 & 1765.328 & 1800 & 1800 & 1800 & 1800 & 1800 \\
 \hline
 Instance 3 & 6.375 & 97.516 & 1800 & 1800 & 1800 & 1800 & 1800 & 1800 \\
 \hline
 Instance 4 & 2.984 & 117.078 & 839.344 & 1800 & 1800 & 1800 & 1800 & 1800 \\
 \hline
 Instance 5 & 3.75 & 221.328 & 1800 & 1800 & 1800 & 1800 & 1800 & 1800 \\
 \hline
\end{tabular}
\end{table}

\begin{table}[h!] \footnotesize
\caption {Optimal gap for Program~\eqref{Program:mip} over multiple periods}
\begin{tabular}{| c | c | c | c | c | c | c | c | c |} 
 \hline
  & \multicolumn{8}{|c|}{Budget} \\
  & 0 & 100 & 200 & 300 & 400 & 500 & 600 & 700\\
 \hline
 \multicolumn{9}{|c|}{100 Points} \\
 \hline
 Instance 1 & 0.00\% & 0.01\% & 0.01\% & 0.01\% & 0.01\% & 0.01\% & 0.00\% &0.00\% \\
 \hline
 Instance 2 & 0.00\% &  0.01\% &  0.01\% &  0.09\% &  0.01\% &  0.01\% & 0.00\% & 0.00\% \\
 \hline
 Instance 3 & 0.00\% & 0.01\% & 0.01\% & 0.01\% & 0.01\% & 0.01\% & 0.00\% & 0.00\% \\
 \hline
 Instance 4 & 0.00\% & 0.01\% & 0.01\% & 0.01\% & 0.01\% & 0.00\%& 0.00\% & 0.00\% \\
 \hline
 Instance 5 & 0.00\% & 0.01\% & 0.01\%  & 0.01\%  & 0.14\%  & 0.05\% & 0.00\% & 0.00\% \\
 \hline
 \multicolumn{9}{|c|}{150 Points} \\
 \hline
 Instance 1 & 0.00\% & 0.01\% & 0.01\% & 0.01\% & 0.95\% & 0.40\% & 0.03\% & 0.00\% \\
 \hline
 Instance 2 & 0.00\%& 0.01\% & 0.01\% & 0.01\% & 0.44\% &0.44\%& 0.07\% & 0.00\% \\
 \hline
 Instance 3 & 0.00\% & 0.01\% &0.01\% & 0.23\%& 1.56\% &1.17\% &0.37\% & 0.03\% \\
 \hline
 Instance 4 &0.00\% & 0.01\% & 0.01\% & 0.01\% &0.66\% & 0.27\% & 0.01\% & 0.00\% \\
 \hline
 Instance 5 & 0.00\% & 0.01\% & 0.09\% & 0.03\% & 0.07\% &0.45\% & 0.13\% &0.01\% \\
 \hline
 \multicolumn{9}{|c|}{200 Points} \\
 \hline
 Instance 1 & 0.00\% &0.01\% &0.31\% & 2.34\% & 2.33\% & 1.90\% &0.66\% & 0.22\% \\
 \hline
 Instance 2 & 0.00\% & 0.01\% & 0.01\% & 2.26\% & 2.24\% & 1.82\% & 0.84\% & 0.24\% \\
 \hline
 Instance 3 & 0.00\% & 0.01\%& 0.34\% &3.21\% &3.69\% & 2.85\% & 1.52\% & 0.54\% \\
 \hline
 Instance 4 & 0.00\% & 0.01\% & 0.01\% & 1.15\% & 0.97\% &0.99\% & 0.33\% & 0.13\% \\
 \hline
 Instance 5 & 0.00\% & 0.01\% & 0.46\%& 2.79\% & 2.73\% & 1.89\% & 0.88\% & 0.28\%\\
 \hline
\end{tabular}\label{table:A12}
\end{table}



\end{document}